# TRANSIENT RANDOM WALKS ON A STRIP IN A RANDOM ENVIRONMENT[1]

By Alexander Roitershtein

*Iowa State University*

We consider transient random walks on a strip in a random environment. The model was introduced by Bolthausen and Goldsheid [*Comm. Math. Phys.* **214** (2000) 429–447]. We derive a strong law of large numbers for the random walks in a general ergodic setup and obtain an annealed central limit theorem in the case of uniformly mixing environments. In addition, we prove that the law of the "environment viewed from the position of the walker" converges to a limiting distribution if the environment is an i.i.d. sequence.

**1. Introduction.** In this paper we consider *random walks in a random environment* (RWRE for short) on the strip $\mathbb{Z} \times \{1, \ldots, d\}$ for some fixed $d \in \mathbb{N}$. Transition probabilities of the random walk are not homogeneous in the space and depend on a realization of the *environment*. The environment is a random sequence $\omega = (\omega_n)_{n \in \mathbb{Z}}$, and, given $\omega$, the random walk $X_n$ is a time-homogeneous Markov chain on the strip with transition kernel $H_\omega$ such that:

(i) $H_\omega((n, i), (m, j)) = 0$ if $|n - m| > 1$, that is, transitions from a site $(n, i)$ are only possible either within the layer $n \times \{1, \ldots, d\}$ or to the two neighbor layers;

(ii) for $n \in \mathbb{Z}$, the $d \times d$ matrices $H_\omega((n, \cdot), (m, \cdot))$ with $m = n - 1$, $n$, $n + 1$ are functions of the random variable $\omega_n$.

One can study properties of the random walk $X_n$ that hold either for almost every realization of the environment (so-called *quenched* setting) or under

---

the law obtained by an averaging over the set of environments (the *annealed* setup).

The RWRE on a strip were introduced by Bolthausen and Goldsheid [5]. It is observed in [5] that the RWRE on $\mathbb{Z}$ with bounded jumps introduced by Key [19] can be viewed as a special case of this model. Other particular cases include a transformation of directed-edge-reinforced random walks on graphs [16], persistent RWRE on $\mathbb{Z}$ [1, 25] and, more generally, certain RWRE on $\mathbb{Z}$ with a finite memory. In the last case, the elements of the set $\{1, \ldots, d\}$ represent states of the memory.

Due to the uniqueness of the escape direction in the state space of the random walk, the model shares some common features with the well studied nearest-neighbor RWRE on $\mathbb{Z}$ (see, e.g., [27], Section 2, for a comprehensive survey). However, due to the inhomogeneity of the kernel $H_\omega$ in the second coordinate, there are important differences that make standard one-dimensional techniques not directly available for the study of the RWRE on strips. For instance, in contrast to the one-dimensional model, for a general RWRE on a strip (1) exit probabilities of the random walk are rather nonexplicit functionals of the environment; (2) excursions of the walker between layers are dependent on each other and do not form any probabilistic branching structure; (3) the random times $\tau_n = T_n - T_{n-1}$, where $T_n$ is the first hitting time of the layer $n$, are not independent in the quenched setting and are not stationary in the annealed setup. Some of these issues are addressed in [5] and in the recent work of Goldsheid [15], as well as in the present paper (see, e.g., Lemma 3.2 below).

For RWRE with bounded jumps on $\mathbb{Z}$ it is possible to transform potential equations associated with the random walk into evolution equations of a linear dynamical system. This approach for establishing recurrence criteria was initiated by Key [19] and further developed by Lëtchikov [21, 22, 23], Derriennic [9] and Brémont [7, 8]; see also the paper of Keane and Rolles [25]. Using a different method, the recurrence and transience behavior of the random walks on a strip is related in [5] to the sign of the top Lyapunov exponent of a sequence of random matrices $a_n$ [defined below by (8)].

The main results of this paper are a strong law of large numbers (cf. Theorem 2.1) and an annealed central limit theorem (cf. Theorem 2.5) for the first coordinate of transient RWRE on a strip. The asymptotic speed of the random walk is expressed in terms of certain matrices (in particular $a_n$) introduced in [5] and we give a sufficient condition for the nonzero speed regime (cf. Corollary 2.2). To prove the law of large numbers and the limit theorem we introduce an annealed measure that makes $\tau_n = T_n - T_{n-1}$ into a stationary ergodic sequence. The results that we obtain are similar to, although are less explicit, than the corresponding statements for nearest-neighbor RWRE on $\mathbb{Z}$. In two particular cases, namely for persistent RWRE on $\mathbb{Z}$ (cf. [1, 25]) and for RWRE on $\mathbb{Z}$ with bounded jumps (cf. [7, 8, 21]), similar results were previously obtained by different methods. We note that a quenched CLT for RWRE on a strip was recently obtained in [15]. In addition, for a class of i.i.d. environments we prove that the law of the "environment viewed from the position of the walker" converges toward a limiting distribution. Similar results for RWRE on $\mathbb{Z}^d$ were previously established by Kesten [17] ($d = 1$) and Sznitman and Zerner [26] ($d > 1$). The renewal structure that we use in order to prove this theorem is a hybrid of the two introduced in [17] and [26].

We turn now to a precise definition of the model. Let $d \geq 1$ be any integer and denote $\mathcal{D} = \{1, \ldots, d\}$. The environment $\omega = (\omega_n)_{n \in \mathbb{Z}}$ is a stationary and ergodic sequence of random variables taking values in a measurable space $(\mathcal{S}, \mathcal{B})$. We denote the distribution of $\omega$ on $(\Omega, \mathcal{F}) := (\mathcal{S}^{\mathbb{Z}}, \mathcal{B}^{\otimes \mathbb{Z}})$ by $P$ and the expectation with respect to $P$ by $E_P$.

Let $\mathcal{A}_d$ be the set of $d \times d$ matrices with real-valued nonnegative entries and let $\mathcal{P}_d \subset \mathcal{A}_d$ denote the set of stochastic matrices. Let three functions $p = (p(i, j))_{i, j \in \mathcal{D}}$, $q = (q(i, j))_{i, j \in \mathcal{D}}$, $r = (r(i, j))_{i, j \in \mathcal{D}} : \mathcal{S} \to \mathcal{A}_d$ be given such that $p + q + r \in \mathcal{P}_d$. For $n \in \mathbb{Z}$, define the following random variables in $(\Omega, \mathcal{F}, P) : p_n = p(\omega_n)$, $q_n = q(\omega_n)$, $r_n = r(\omega_n)$. The triples $(p_n, r_n, q_n)$ will facilitate the definition of transition probabilities of the random walk $X_n$.

We will denote the first component of $X_n$ by $\xi_n$ and the second one by $Y_n$, that is,

$$X_n = (\xi_n, Y_n), \qquad \xi_n \in \mathbb{Z}, \ Y_n \in \mathcal{D}.$$

It will be convenient to consider random walks for which $\xi_0 = 0$ with probability 1, but the initial position $Y_0$ of the second component is random and its distribution is a function of the environment $\omega_-(\omega) := (\omega_n)_{n \leq 0}$. Let $\mathcal{P}r_d \subset \mathbb{R}_+^d$ be the set of probability measures on the finite set $\mathcal{D}$ and define $\mathcal{M}_d \subset (\mathcal{P}r_d)^\Omega$ as follows:

$$\mathcal{M}_d = \{(\mu_\omega)_{\omega \in \Omega} : \omega \to \mu_\omega \text{ is a measurable function from } \Omega \text{ to } \mathcal{P}r_d$$

$$\text{and } \mu_\alpha = \mu_\beta \text{ if } \omega_-(\alpha) = \omega_-(\beta)\}.$$

Given a collection $\mu = (\mu_\omega)_{\omega \in \Omega} \in \mathcal{M}_d$, of initial distributions on $\mathcal{D}$, the *random walk on the strip* $\mathbb{Z} \times \mathcal{D}$ *in environment* $\omega \in \Omega$ is a time-homogeneous Markov chain $X = (X_t)_{t \in \mathbb{Z}_+}$ taking values in $\mathbb{Z} \times \mathcal{D}$ and governed by the *quenched* law $P_\omega^\mu$ defined by transition probabilities (here $n \in \mathbb{Z}$ and $i, j \in \mathcal{D}$)

$$H_\omega(x, y) := \begin{cases} p_n(i, j), & \text{if } x = (n, i) \text{ and } y = (n + 1, j), \\ r_n(i, j), & \text{if } x = (n, i) \text{ and } y = (n, j), \\ q_n(i, j), & \text{if } x = (n, i) \text{ and } y = (n - 1, j), \\ 0, & \text{otherwise,} \end{cases}$$

and initial distribution $P_\omega^\mu(\xi_0 = 0, Y_0 = y_0) = \mu_\omega(y_0)$ for any $y_0 \in \mathcal{D}$. That is, for any finite sequence $x_0 = (0, y_0), x_1, \ldots, x_n$ of elements in $\mathbb{Z} \times \mathcal{D}$,

$$P_\omega^\mu(X_0 = x_0, X_1 = x_1, \ldots, X_n = x_n)$$

$$= \mu_\omega(y_0) H_\omega(x_0, x_1) \cdots H_\omega(x_{n-1}, x_n).$$

Let $\mathcal{G}$ be the cylinder $\sigma$-algebra on $(\mathbb{Z} \times \mathcal{D})^{\mathbb{N}}$, the path space of the walk. The *random walk (associated with $\mu$) on $\mathbb{Z} \times \mathcal{D}$ in a random environment* is the process $(\omega, X)$ in the measurable space $(\Omega \times (\mathbb{Z} \times \mathcal{D})^{\mathbb{N}}, \mathcal{F} \times \mathcal{G})$ with the *annealed* law $\mathbb{P}^{\mu} = P \otimes P_{\omega}^{\mu}$ defined by

$$\mathbb{P}^{\mu}(F \times G) = \int_F P_{\omega}^{\mu}(G) P(d\omega) = E_P(P_{\omega}^{\mu}(G); F), \qquad F \in \mathcal{F}, \ G \in \mathcal{G}.$$

That is, the annealed law of $X_n$ is obtained by averaging its quenched law (i.e., given a fixed environment) over the set of environments.

Let us introduce some further notation. We denote by $P_{\omega}^i$ (resp., $\mathbb{P}^i$) the quenched (annealed) law of the random walk starting at site $(0, i)$. That is, $\mathbb{P}^i(\cdot) = E_P(P_{\omega}^i(\cdot))$ and

$$P_{\omega}^i(\cdot) = P_{\omega}^{\mu}(\cdot) \qquad \text{with } \mu_{\omega}(i) = 1.$$

Throughout the paper we use the notation $\mathbf{0} := (0, \ldots, 0) \in \mathbb{R}^d$, $\mathbf{1} := (1, \ldots, 1) \in \mathbb{R}^d$, and denote by $e_i$, $i = 1, \ldots, d$, the vectors $(0, \ldots, 1, \ldots, 0)$ of the canonical basis of $\mathbb{R}^d$. We denote by $O$ the zero $d \times d$ matrix and by $I$ the unit $d \times d$ matrix. The notions ">," "<," "positive," "negative" and "nonnegative" are used for vectors and matrices in the usual way. For instance, a $d \times d$ matrix $A$ is said to be nonnegative if $A \geq O$, that is, $A(i, j) \geq 0$ for all $i, j \in \mathcal{D}$. For a vector $x = (x^1, \ldots, x^d) \in \mathbb{R}^d$ let $\|x\| := \max_{i \in \mathcal{D}} |x^i|$ and for a $d \times d$ real matrix $A$ let $\|A\|$ denote the corresponding operator norm $\sup_{\{x \in \mathbb{R}^d : \|x\| = 1\}} \|Ax\| = \max_{i \in \mathcal{D}} \sum_{j \in \mathcal{D}} |A_{i,j}|$. We frequently use the fact that if a $d \times d$ matrix $A$ is nonnegative, then $\|A\| = \|A\mathbf{1}\|$. Finally, we use the notation $\mathbf{I}_A$ for the indicator function of the set $A$.

We next introduce our basic assumptions. For $n \in \mathbb{Z}$, define the random times

$$(1) \qquad T_n = \inf\{i : \xi_i = n\} \quad \text{and} \quad \tau_n := T_n - T_{n-1},$$

with the usual convention that the infimum over an empty set is $\infty$ and $\infty - \infty = \infty$. Let $\eta_n(i, j)$ be the probability that the walker starting at site $(n, i)$ in environment $\omega$ reaches the layer $n + 1$ at point $(n + 1, j)$, that is (with some abuse of notation, in this and similar cases, we usually will not indicate explicitly the dependence on $\omega$),

$$(2) \qquad \eta_n(i, j) = P_{\theta^n \omega}^i(Y_{\tau_1} = j),$$

where $\theta : \Omega \to \Omega$ is the shift operator defined by

$$(3) \qquad (\theta \omega)_n = \omega_{n+1}, \qquad n \in \mathbb{Z}.$$

A "continued fractions" representation of $\eta_n$ in terms of the sequence $\omega_n$ is provided in [5], page 437 [see especially identities (2.5) and (2.7) there].

The following Condition C is borrowed from the paper of Bolthausen and Goldsheid [5].

CONDITION C.

(C1) The sequence $(\omega_n)_{n \in \mathbb{Z}}$ is stationary and ergodic.

(C2) $E_P(\log(1 - \|r_0 + p_0\|)^{-1}) < \infty$, $E_P(\log(1 - \|r_0 + q_0\|)^{-1}) < \infty$.

(C3) For all $j \in \mathcal{D}$,

$$\sum_{i=1}^{d} q_0(i,j) > 0 \quad \text{and} \quad \sum_{i=1}^{d} p_0(i,j) > 0,$$

$P$-almost surely.

(C4) $P(\eta_0(i,j) > 0) = 1$ for all $i, j \in \mathcal{D}$.

Note that since

$$(4) \qquad \min_{i \in \mathcal{D}} \sum_{j \in \mathcal{D}} q_0(i,j) = 1 - \|r_0 + p_0\|, \qquad \min_{i \in \mathcal{D}} \sum_{j \in \mathcal{D}} p_0(i,j) = 1 - \|r_0 + q_0\|.$$

Condition C2 implies that for all $i \in \mathcal{D}$, $P$-a.s.,

$$(5) \qquad \sum_{j=1}^{d} q_0(i,j) > 0 \quad \text{and} \quad \sum_{j=1}^{d} p_0(i,j) > 0.$$

In this paper we concentrate on random walks $X_n$ which are *transient to the right*, that is, $\mathbb{P}^i(\lim_{n \to \infty} \xi_n = \infty) = 1$ for all $i \in \mathcal{D}$. We next recall the criterion for the transient behavior of $X_n$ obtained by Bolthausen and Goldsheid in [5]. Let $\zeta_n = \zeta_n(\omega)$, $n \in \mathbb{Z}$, be a stationary sequence of $d \times d$ stochastic matrices solving the equation

$$\zeta_n = (I - q_n \zeta_{n-1} - r_n)^{-1} p_n, \qquad P\text{-a.s.}, \qquad n \in \mathbb{Z}.$$

This sequence exists and is in fact unique by Theorem 1 in [5]. Let $\gamma_n = (I - q_n \zeta_{n-1} - r_n)^{-1}$, and introduce for $n \in \mathbb{Z}$ the following random matrices: $a_n = \gamma_n q_n = (I - q_n \zeta_{n-1} - r_n)^{-1} q_n$. As we shall see, the matrices $a_n$ play essentially the same role for RWRE on a strip as the random variables $\rho_n = P_\omega(\text{jump from } n \text{ to } n-1)/P_\omega(\text{jump from } n \text{ to } n+1)$ do for the nearest-neighbor RWRE on $\mathbb{Z}$.

By the subadditive ergodic theorem of [14], condition C2 ensures that there exists a nonrandom number $\lambda$ such that

$$(6) \qquad \begin{aligned} \lambda &= \lim_{n \to \infty} 1/n \cdot E_P(\log \|a_n \cdots a_1\|) \\ &= \lim_{n \to \infty} 1/n \cdot \log \|a_n \cdots a_1\|, \qquad P\text{-a.s.} \end{aligned}$$

It is shown in [5] that under Condition C, $X_n$ is transient to the right if and only if $\lambda < 0$. Therefore we will assume throughout the paper:

CONDITION D. $\lambda < 0$.

This condition is a generalization of the transience criterion $E_P(\log \rho_0) < \infty$ for the RWRE on $\mathbb{Z}$ (see, e.g., [27], Section 2). If the random walk is transient to the right, $P(\eta_n = \zeta_n) = 1$ [5] and hence we have $P$-a.s. for all $n \in \mathbb{Z}$:

$$\gamma_n = (I - q_n \eta_{n-1} - r_n)^{-1} \tag{7}$$

and

$$a_n = (I - q_n \eta_{n-1} - r_n)^{-1} q_n. \tag{8}$$

The rest of the paper is organized as follows. The main results are stated in Section 2. A law of large numbers for $\xi_n$ (cf. Theorem 2.1) and a criterion for a positive asymptotic speed (cf. Corollary 2.2) are proved in Section 3. Section 4 is devoted to the Markov chain "environment viewed from the particle." In particular, we prove there that the law of the Markov chain converges toward a limiting distribution. An annealed central limit theorem for a normalized random variable $\xi_n$ in the case of uniformly mixing environments is given by Theorem 2.5 and is proved in Section 5.

**2. Statement of main results.** Let $\pi = (\pi_\omega)_{\omega \in \Omega} \in \mathcal{M}_d$ be a collection of probability measures on $\mathcal{D}$ defined $P$-a.s. by

$$\pi_\omega = \lim_{n \to \infty} e_i \eta_{-n}(\omega) \cdots \eta_{-1}(\omega). \tag{9}$$

It follows, for example, from part (iii) of Lemma 3.3 that the limit in (9) exists $P$-a.s., does not depend on $i$, and defines a probability measure on $\mathcal{D} = \{1, \ldots, d\}$ whose support is the whole $\mathcal{D}$.

It is shown in Section 3 that the sequence $(\tau_n)_{n \in \mathbb{N}}$ is stationary and ergodic under $\mathbb{P}^\pi$ and the following law of large numbers holds.

THEOREM 2.1. *Let Conditions C and D hold. Then for every $\mu \in \mathcal{M}_d$, the following limits exist $\mathbb{P}^\mu$-a.s. and the equality holds:*

$$v_P := \lim_{n \to \infty} \xi_n/n = \lim_{n \to \infty} n/T_n$$
$$= 1/E_P(\pi_\omega \cdot (b_0 + a_0 b_{-1} + a_0 a_{-1} b_{-2} + \cdots)),$$

*where $b_n := \gamma_n \mathbf{1}$.*

The formula for the asymptotic speed yields the following characterization of the positive speed regime.

COROLLARY 2.2. *Assume that:*

(i) *Condition C holds.*

(ii) *There exists $\epsilon > 0$ such that $P$-a.s., $\sum_{j=1}^d p_n(i,j) > \epsilon$ for all $i \in \mathcal{D}$ and all $n \in \mathbb{Z}$.*

(iii) *$\limsup_{n \to \infty} 1/n \log E_P(\|a_0 a_1 \cdots a_{n-1}\|) < 0$.*

*Then $P(v_P > 0) = 1$.*

Assumption (iii) of Corollary 2.2 is a generalization of the corresponding condition $\limsup_{n \to \infty} \frac{1}{n} \log E_P(\rho_0 \rho_1 \cdots \rho_{n-1}) < 0$ for RWRE on $\mathbb{Z}$ (cf. [27], Section 2). Note that by Jensen's inequality this assumption implies Condition D. Assumption (ii) implies in virtue of (4) that $P$-a.s.,

$$
\begin{aligned}
\|b_n\| = \|\gamma_n \mathbf{1}\| = \|\gamma_n\| &\leq (1 - \|q_n \eta_{n-1} + r_n\|)^{-1} \\
&= (1 - \|q_n + r_n\|)^{-1} \leq 1/\epsilon,
\end{aligned}
\tag{10}
$$

and hence, by condition (iii) of the corollary,

$$
\begin{aligned}
\limsup_{n \to \infty} \frac{1}{n} &\log E_P(\|a_0 a_1 \cdots a_{-n+1} b_{-n}\|) \\
&\leq \limsup_{n \to \infty} 1/n \log E_P(\|a_0 a_1 \cdots a_{-n+1}\| \cdot \|b_{-n}\|) < 0.
\end{aligned}
\tag{11}
$$

Consequently, $v_P > 0$ because the formula for the asymptotic speed given by Theorem 2.1 and inequality (11) imply that

$$
\begin{aligned}
\frac{1}{v_P} &\leq \sum_{n=0}^{\infty} E_P(\|\pi_\omega\| \cdot \|a_0 a_1 \cdots a_{-n+1} b_{-n}\|) \\
&\leq \sum_{n=0}^{\infty} E_P(\|a_0 a_1 \cdots a_{-n+1} b_{-n}\|) \leq \sum_{n=0}^{\infty} a e^{-bn} < \infty
\end{aligned}
$$

for some constants $a, b > 0$.

In Section 4 we study the "environment viewed from the particle" process $(\theta^{\xi_n}\omega, Y_n)_{n \geq 0}$. By Lemmas 3.1 and 3.2, $v_P = 1/\mathbb{E}^\pi(T_1)$. Therefore Propositions 4.1 and 4.8 yield

THEOREM 2.3. *Let Conditions C and D hold. Then $v_P > 0$ if and only if there is a stationary distribution $Q$ for the Markov chain $x_n = (\theta^{\xi_n}\omega, Y_n)$, $n \geq 0$, which is equivalent to $P \otimes \{$counting measure on $\mathcal{D}\}$. If such a distribution $Q$ exists it is unique and is given by the following formula:*

$$
Q(B, i) = v_P \mathbb{E}^\pi \left( \sum_{n=0}^{T_1 - 1} \mathbf{I}_B(\overline{\omega}_n) \mathbf{I}_{Y_n}(i) \right), \qquad B \in \mathcal{F}, \ i \in \mathcal{D},
\tag{12}
$$

*where $\overline{\omega}_n := \theta^{\xi_n}\omega$.*

For i.i.d. environments, similarly to the RWRE on $\mathbb{Z}^d$ (see [17] for $d = 1$ and [26] for $d > 1$), the law of the environment viewed from the particle converges toward its unique stationary distribution. More precisely, we prove in Section 4 the following theorem.

THEOREM 2.4.   *Assume that $(\omega_n)_{n \in \mathbb{Z}}$ is an i.i.d. sequence, Conditions $C$ and $D$ hold, and $v_P > 0$. Then for any $\mu \in \mathcal{M}_d$, under the law $\mathbb{P}^\mu$ the sequence $x_n = (\theta^{\xi_n}\omega, Y_n)$ converges weakly to the limiting distribution $Q$ as $n$ goes to infinity.*

In the case where $(\omega_n)_{n \in \mathbb{Z}}$ is a uniformly mixing sequence we obtain in Section 5 the following central limit theorem. Let $\sigma^{(n)} = \sigma(\omega_i : i \geq n)$, $\sigma_n = \sigma(\omega_i : i \leq n)$, and

$$\phi(n) = \sup\{P(A|B) - P(A) : A \in \sigma^{(n)}, B \in \sigma_0, P(B) > 0\}.$$

Recall that the sequence $\omega_n$ is called *uniformly mixing* if $\phi(n) \to_{n \to \infty} 0$. Let $\varphi(n) := 2\phi(n)$. For uniformly mixing sequence it holds that (cf. [11], page 9)

$$(13) \qquad |E_P(fg) - E_P(f)E_P(g)| \leq \varphi(n)\|f\|_\infty \|g\|_1$$

for bounded $\sigma(\omega_i : i \geq n)$-measurable functions $f$ and bounded $\sigma(\omega_i : i \leq 0)$-measurable functions $g$.

THEOREM 2.5.   *Let Condition $C$ hold and assume in addition that:*

   (i)  *There exists $\epsilon > 0$ such that $P$-a.s., $\sum_{j=1}^d p_n(i,j) > \epsilon$ for all $i \in \mathcal{D}$ and all $n \in \mathbb{Z}$.*

   (ii)  $\limsup_{n \to \infty} n^{-1} \log E_P(\|a_0 a_1 \cdots a_{n-1}\|^2) < 0$.

   (iii)  $\sum_{n=1}^\infty \sqrt{\varphi(n)} < \infty$, *where the mixing coefficients $\varphi(n)$ are defined in* (13).

   (iv)  *We have $\limsup_{n \to \infty} n^{-1} \log E_P(c_0 c_1 \cdots c_{n-1}) < 0$ for $c_n = 1 - \max_{i \in \mathcal{D}} \min_{j \in \mathcal{D}} \eta_n(i,j)$.*

*Then there exists a positive constant $\sigma > 0$ such that for any $\mu \in \mathcal{M}_d$, $n^{-1/2}(\xi_n - n \cdot v_P) \Rightarrow N(0, \sigma^2)$ under $\mathbb{P}^\mu$ as $n \to \infty$, where $N(0, \sigma^2)$ stands for the normal distribution with variance $\sigma^2$.*

Assumption (ii) of the theorem implies by Jensen's inequality Condition D. Moreover, by Corollary 2.2, $v_P > 0$. We note that condition (iv) of Theorem 2.5 is implied by a number of more explicit assumptions, for instance by the following one: $\limsup_{n \to \infty} n^{-1} \log E_P(\tilde{c}_0 \cdots \tilde{c}_1 \cdots \tilde{c}_{n-1}) < 0$, where $\tilde{c}_n = 1 - \min_{j \in \mathcal{D}} p_n(i^*, j)$ for some $i^* \in \mathcal{D}$.

**3. Law of large numbers for $\xi_n$.** This section is devoted to the proof of Theorem 2.1. The law of large numbers for $\xi_n$ is derived from the corresponding statement for the hitting times $T_n$ introduced in (1). More precisely, we have:

LEMMA 3.1. *Let Conditions C and D hold and assume in addition that* $\mathbb{P}^\mu(\lim_{n\to\infty} T_n/n = \alpha) = 1$ *for some* $\alpha \in (1,\infty]$ *and* $\mu \in \mathcal{M}_d$. *Then* $\mathbb{P}^\mu(\lim_{n\to\infty} \xi_n/n = 1/\alpha) = 1$.

The proof of the lemma is similar to that of the analogous claim for the RWRE on $\mathbb{Z}$ (cf. [27], Lemma 2.1.17) and thus is omitted.

In contrast to the nearest-neighbor RWRE on $\mathbb{Z}$, the sequence $\tau := (\tau_n)_{n\in\mathbb{N}}$ defined in (1) is in general not stationary for transient random walks on the strip. However, as the next lemma shows, the sequence is stationary and ergodic under $\mathbb{P}^\pi$, where $\pi = (\pi_\omega)_{\omega\in\Omega}$ is defined in (9).

LEMMA 3.2. (a) $\tau = (\tau_n)_{n\in\mathbb{N}}$ *is a stationary sequence under* $\mathbb{P}^\pi$.
(b) $\tau = (\tau_n)_{n\in\mathbb{N}}$ *is ergodic under* $\mathbb{P}^\pi$.

Before we turn to the proof of the lemma, let us have a closer look at the matrices $\eta_n$ introduced in (2). Conditioning on the direction of the first step of the random walk we get that $\eta_n$ verify the following equation:

$$\eta_n = p_n + r_n\eta_n + q_n\eta_{n-1}\eta_n,$$

which is equivalent to

$$(14) \qquad \eta_n = (I - q_n\eta_{n-1} - r_n)^{-1} p_n = \gamma_n p_n, \qquad n \in \mathbb{Z}.$$

In fact, under Conditions C and D, $(\eta_n)_{n\in\mathbb{Z}}$ is the unique stationary sequence which satisfies recursion (14) (cf. [5]). In the following lemma we collect some properties of the matrices $\eta_n$ that were proved in [5].

LEMMA 3.3 [5]. *Let Condition C hold and assume in addition that* $\mathbb{P}(\lim_{n\to\infty} \xi_n = +\infty) = 1$. *Then:*

(i) *The extended sequence* $(\omega_n, \eta_n)_{n\in\mathbb{Z}}$ *is stationary and ergodic.*

(ii) *The matrices $\eta_n$ are $P$-a.s. stochastic and $P(\eta_0 > 0) = 1$.*

(iii) *There are an ergodic stationary sequence of matrices $(\pi_n(\omega))_{n\in\mathbb{Z}}$ and a sequence of matrices $(\varepsilon_{n,m}(\omega))_{n\in\mathbb{Z},m\in\mathbb{N}}$ such that for all $n \in \mathbb{Z}$ and $m \in \mathbb{N}$:*

  (a) $\eta_{n-m}\cdots\eta_{n-2}\eta_{n-1} = \pi_n + \varepsilon_{n,m}$.

  (b) *$\pi_n$ are strictly positive column matrices, that is, $P$-a.s., $\pi_n(i,k) = \pi_n(j,k) > 0$ for any $i,j,k \in \mathcal{D}$ and $n \in \mathbb{Z}$.*

  (c) *$P(\|\varepsilon_{n,m}\| \le c_{n-1}\cdots c_{n-m}) = 1$, where $c_n$ are introduced in the statement of Theorem 2.5.*

Part (i) of the above lemma follows from [5], Theorem 1c, because if $\mathbb{P}(\lim_{n \to \infty} \xi_n = +\infty) = 1$, then $\eta_n$ coincides with $\zeta_n$ defined there. Part (ii) is the content of [5], Lemma 4 and [5], Corollary 2. Part (iii) is essentially [5], Lemma 9 with the only exception that our choice of the numbers $c_n$ is different from that of [5]. The above version follows from a contraction property of Doeblin's stochastic kernels (in particular, strictly positive stochastic matrices); compare [10], page 197.

Proof of Lemma 3.2. (a) Recall the shift operator $\theta : \Omega \to \Omega$ defined in (3). By the definition, $\pi_{\theta^n \omega} = \pi_{\theta^{n-1} \omega} \eta_{n-1}$ and hence, by induction on $n$, $P$-a.s. we have $\pi_{\theta^n \omega}(\cdot) = P_\omega^\pi(\xi_{T_n} \in \cdot)$ for $n \geq 0$. Since $(\eta_n(\omega))_{n \in \mathbb{Z}}$ is a stationary sequence under $P$, so is the sequence $(\pi_{\theta^n \omega})_{n \in \mathbb{Z}}$. Let $\theta_\tau$ be the shift operator for the sequence $\tau$ defined by $(\theta_\tau \tau)_n = \tau_{n+1}$, that is, $\theta_\tau \tau(\omega) = \tau(\theta \omega)$. Then, for any $n \in \mathbb{N}$,

$$\mathbb{P}^\pi(\theta_\tau^n \tau \in \cdot) = \sum_{i \in \mathcal{D}} E_P(P_\omega^\pi(Y_{T_n} = i) P_{\theta^n \omega}^i(\tau \in \cdot))$$

$$= \sum_{i \in \mathcal{D}} E_P(\pi_{\theta^n \omega}(i) P_{\theta^n \omega}^i(\tau \in \cdot))$$

$$= \sum_{i \in \mathcal{D}} E_P(\pi_\omega(i) P_\omega^i(\tau \in \cdot)) = \mathbb{P}^\pi(\tau \in \cdot),$$

where the last but one step follows from the translation invariance of the measure $P$ with respect to the shift $\theta$.

(b) For $n \in \mathbb{N}$, let $\bar{x}_n = (\theta^n \omega, Y_{T_n}, \tau_n)$. It is easy to see that under $\mathbb{P}^\pi$, $(\bar{x}_n)_{n \in \mathbb{N}}$ is a stationary Markov chain with transition kernel

$$K(\omega, i, n; A, j, m) = P_\omega^i(\tau_1 = m, Y_m = j) \mathbf{I}_A(\theta \omega).$$

The claim of part (b) follows from the following lemma:

Lemma 3.4. *The sequence $(\bar{x}_n)_{n \in \mathbb{N}}$ is ergodic under $\mathbb{P}^\pi$.*

We turn now to the proof of Lemma 3.4. Denote $\mathcal{X} := (\bar{x}_n)_{n \in \mathbb{Z}_+}$ and let $\mathbb{S} = \Omega \times \mathcal{D} \times \mathbb{N}$ be the state space of the Markov chain $(\bar{x}_n)_{n \in \mathbb{N}}$. Let $A \in \sigma(\bar{x}_n : n \geq 0)$ be an invariant set, that is, $\eta^{-1} A = A$, $\mathbb{P}^\pi$-a.s., where $\eta$ is the usual shift in the space $\mathbb{S}^{\mathbb{Z}_+} : (\eta \mathcal{X})_n = \bar{x}_{n+1}$. It suffices to show that $\mathbb{P}^\pi(A) \in \{0, 1\}$. For $x \in \mathbb{S}$ let $\mathbb{P}_x$ be the law of the Markov chain $(\bar{x}_n)_{n \geq 0}$ with the initial state $\bar{x}_0 = x$. Set $h(x) = \mathbb{P}_x(A)$. The following lemma is a general property of stationary Markov chains (see, e.g., [6], page 12 or [27], page 207):

Lemma 3.5. (i) *The sequence $(h(\bar{x}_n))_{n \in \mathbb{N}}$ forms a martingale in its canonical filtration.*

(ii) *There exists a measurable set $B \subset \mathbb{S}$ such that $h(\bar{x}_0) = \mathbf{I}_B(\bar{x}_0)$, $\mathbb{P}^\pi$-a.s.*

It remains to show that $\mathbb{P}^\pi(\bar{x}_0 \in B) \in \{0,1\}$. By the martingale property of the sequence $h(\bar{x}_n)$, $\mathbf{I}_B(x) = K\mathbf{I}_B(x)$, we have $\mathbb{P}^\pi$-a.s. Write $B = \bigcup_{i\in\mathcal{D}}\bigcup_{n\in\mathbb{N}} B_{i,n} \times \{i\} \times \{n\}$. We have, for $P$-a.e. environment $\omega$,

$$\mathbf{I}_{B_{i,n}}(\omega) = \mathbf{I}_B(\omega, i, n) = K\mathbf{I}_B(\omega, i, n)$$
$$= \sum_{j\in\mathcal{D}}\sum_{m\in\mathbb{N}} P_\omega^i(\tau_1 = m, Y_m = j)\mathbf{I}_{B_{j,m}}(\theta\omega).$$

Therefore, for all $m, n \in \mathbb{N}$, $\mathbf{I}_{B_{i,n}}(\omega) = \mathbf{I}_{B_{i,m}}(\omega)$ with probability 1. That is, using the notation $\mathbf{I}_{B_i} := \mathbf{I}_{B_{i,n}}$, $n \in \mathbb{N}$, we can write

$$(15) \qquad B = \bigcup_{i\in\mathcal{D}} B_i \times \{i\},$$

and the preceding equation becomes

$$\mathbf{I}_{B_i}(\omega) = \sum_{j\in\mathcal{D}} P_\omega^i(Y_{T_1} = j)\mathbf{I}_{B_j}(\theta\omega) = \sum_{j\in\mathcal{D}} \eta_0(i,j)\mathbf{I}_{B_j}(\theta\omega).$$

Since $\mathbf{I}_{B_i} \in \{0,1\}$, and by condition C4, $P(\eta_0 > 0) = 1$, we obtain that

$$(16) \qquad \mathbf{I}_{B_i}(\omega) = \mathbf{I}_{B_j}(\theta\omega) \qquad \text{for all } i, j \in \mathcal{D}, \qquad P\text{-a.s.}$$

In particular, $\mathbf{I}_{B_i}(\omega)$ are invariant functions of the ergodic sequence $(\omega_n)_{n\in\mathbb{Z}}$, and hence $P(B_i) \in \{0,1\}$. Furthermore, by (16) and the translation invariance of the measure $P$, we have $P(B_i) = P(\theta B_j) = P(B_j)$ for every $i, j \in \mathcal{D}$. It follows from (15) that $\mathbb{P}^\pi(B) \in \{0,1\}$. The proof of the lemma is complete. $\square$

Applying the ergodic theorem to the sequence $\tau = (\tau_n)_{n\geq 1}$, we obtain that $\mathbb{P}^\pi(\lim_{n\to\infty} T_n/n = \alpha) = 1$, where $\alpha := \mathbb{E}^\pi(T_1)$. Since $P(\pi_\omega > 0) = 1$ (cf. Lemma 3.3), the a.s. convergence of $T_n/n$ under $\mathbb{P}^\pi$ implies the a.s. convergence under $\mathbb{P}^\mu$ for any $\mu \in \mathcal{M}_d$. In the next lemma we compute the limit $\alpha = 1/v_P$. For any integer $n \leq 0$ and $j \in \mathcal{D}$ let

$$(17) \qquad u_n^j(\omega) = E_\omega(T_1 | X_0 = (n,j)),$$

and denote by $u_n$ the vector $(u_n^1, \ldots, u_n^d)$. Then $\alpha = E_P(\pi_\omega u_0(\omega))$ and we have:

LEMMA 3.6. *Let the conditions of Theorem 2.1 hold. Then $P$-a.s.,*

$$u_0 = b_0 + a_0 b_{-1} + a_0 a_{-1} b_{-2} + \cdots < \infty,$$

*where $b_n := \gamma_n \mathbf{1}$ and the random matrices $\gamma_n$ and $a_n$ are defined in (7) and (8), respectively.*

PROOF. For $n \leq 0$ and $i \in \mathcal{D}$, decompose the path of the random walk with $X_0 = (n, i)$ according to the value of $X_1$:

$$\tau_1 = \sum_{k \in \{-1, 0, 1\}} \sum_{j \in \mathcal{D}} \mathbf{I}_{\{\xi_1 = n+k, Y_1 = j\}} (1 + \tilde{\tau}_1), \tag{18}$$

where $(1 + \tilde{\tau}_1)$ is the (possibly infinite) first hitting time of the layer $\{1\} \times \mathcal{D}$ after time 1.

Setting $u_1 = \mathbf{0}$ and taking expectations in the previous identity, we obtain

$$u_n = \mathbf{1} + q_n u_{n-1} + r_n u_n + p_n u_{n+1}, \qquad n \leq 0. \tag{19}$$

Let $y_1 := u_0$ and, using the usual convention that $\infty - \infty = \infty$, define for $n \leq 0$,

$$y_n = u_{n-1} - \eta_{n-1} u_n. \tag{20}$$

Note that by the definition $y_n \geq 0$ for all $n \leq 0$. Indeed,

$$
\begin{aligned}
u_{n-1}^i &= E_\omega(T_1 | X_0 = (n-1, i)) \\
&= \sum_{j \in \mathcal{D}} \eta_{n-1}(i, j) E_\omega(T_n + T_1 - T_n | X_0 = (n-1, i), Y_{T_n} = j) \\
&\geq \sum_{j \in \mathcal{D}} \eta_{n-1}(i, j) E_\omega(T_1 - T_n | X_0 = (n-1, i), Y_{T_n} = j) \\
&= \sum_{j \in \mathcal{D}} \eta_{n-1}(i, j) E_\omega(T_1 | X_0 = (n, j)) = (\eta_{n-1} u_n)^i,
\end{aligned}
$$

where $(\eta_{n-1} u_n)^i$ denotes the $i$th component of the vector $\eta_{n-1} u_n$.

Substituting the expression for $u_{n-1}$ given by (20) into (19) we obtain

$$(I - q_n \eta_{n-1} - r_n) u_n = \mathbf{1} + q_n y_n + p_n u_{n+1}, \qquad n \leq 0,$$

and hence, in view of (14) and (8),

$$y_{n+1} = b_n + a_n y_n, \qquad n \leq 0. \tag{21}$$

Iterating the identity (21) and using a standard truncation argument, we will next show that

$$u_0 = y_1 = b_0 + a_0 b_{-1} + a_0 a_{-1} b_{-2} + \cdots. \tag{22}$$

First, observe that (21) yields for any $n \in \mathbb{N}$,

$$
\begin{aligned}
y_1 &= b_0 + a_0 b_{-1} + a_0 a_{-1} b_{-2} + \cdots + a_0 \cdots a_{-n+1} b_{-n} + a_0 \cdots a_{-n} y_{-n} \\
&\geq b_0 + a_0 b_{-1} + a_0 a_{-1} b_{-2} + \cdots + a_0 a_{-1} \cdots a_{-n+1} b_{-n},
\end{aligned}
$$

and hence $y_1 \geq b_0 + \sum_{n=1}^{\infty} a_0 a_{-1} \cdots a_{-n+1} b_{-n}$.

To show that the reverse inequality holds, fix any real number $M > 0$ and define for integers $n \leq 0$ $d$-vectors $v_{n,M}$ by setting

$$v_{n,M}^j(\omega) = E_\omega(T_1 \wedge M | X_0 = (n,j))$$

for the $j$th component of the vector. Here we use the standard notation $x \wedge y = \min\{x, y\}$. It follows from (18) that

$$(23) \qquad v_{n,M} \leq \mathbf{1} + q_n v_{n-1,M} + r_n v_{n,M} + p_n v_{n+1,M}, \qquad n \leq 0.$$

Indeed, if at least one of the summands $\mathbf{I}_{\{\xi_1 = \cdot, Y_1 = \cdot\}}(1 + \tilde{\tau}_1)$ in the right-hand side of (18) is bigger than or equal to $M$, the left-hand side is at least $M$ as well. In this case

$$(24) \qquad \begin{aligned} \tau_1 \wedge M &\leq \sum_{k \in \{-1,0,1\}} \sum_{j \in \mathcal{D}} \mathbf{I}_{\{\xi_1 = n+k, Y_1 = j\}} (1 + \tilde{\tau}_1) \wedge M \\ &\leq \mathbf{1} + \sum_{k \in \{-1,0,1\}} \sum_{j \in \mathcal{D}} \mathbf{I}_{\{\xi_1 = n+k, Y_1 = j\}} \tilde{\tau}_1 \wedge M. \end{aligned}$$

Taking the expectations in (24) we obtain inequality (23).

Letting $z_{1,M} := v_{0,M}$ and for $n \leq 0$, $z_{n,M} := v_{n-1,M} - \eta_{n-1} v_{n,M}$, we get from (23), similarly to (21), that

$$z_{n+1,M} \leq b_n + a_n z_{n,M}, \qquad n \leq 0.$$

Since $\|z_{n,M}\| \leq M$ and $\|a_0 \cdots a_{-n+1}\| \to 0$ by Condition D, it follows that $z_{1,M} \leq b_0 + \sum_{n=1}^\infty a_0 a_{-1} \cdots a_{-n+1} b_{-n}$ and hence (22) holds.

To complete the proof of the lemma it remains to show that the right-hand side of (22) is $P$-a.s. finite. In fact, the following standard argument implies that, under the conditions of Theorem 2.1, the norm of the vectors $a_0 a_{-1} \cdots a_{-n} b_{-n-1}$ decay exponentially. Namely,

$$\lim_{n \to \infty} 1/n \cdot \log \|a_0 a_{-1} \cdots a_{-n} b_{-n-1}\| < 0, \qquad P\text{-a.s.}$$

In virtue of Condition D and the multiplicative property of the norm, it suffices to show that

$$(25) \qquad \lim_{n \to \infty} 1/n \cdot \log \|b_n\| = 0, \qquad P\text{-a.s.}$$

Toward this end observe first that

$$\begin{aligned} E_P(\log \|b_n\|) &= E_P(\log \|\gamma_n \mathbf{1}\|) = E_P(\log \|\gamma_n\|) \\ &\leq E_P(\log(1 - \|q_n \eta_{n-1} + r_n\|)^{-1}) \\ &= E_P((1 - \|q_n + r_n\|)^{-1}) < \infty, \end{aligned}$$

where the last inequality is due to Condition C2. Since the sequence of matrices $(\gamma_n)_{n\in\mathbb{Z}}$ defined in (7) is stationary, and hence the sequence of vectors $b_n = \gamma_n \mathbf{1}$ is stationary as well, we have for any fixed $\varepsilon > 0$,

$$\sum_{n=1}^{\infty} P(1/n \log \|b_n\| > \varepsilon) = \sum_{n=1}^{\infty} P(\log \|b_n\| > n\varepsilon)$$
$$= E_P(\log \|b_0\|) - 1 < \infty,$$

implying by the Borel–Cantelli lemma that $P(1/n \log \|b_n\| > \varepsilon \text{ i.o.}) = 0$. □

**4. Environment viewed from the particle.** This section is devoted to the proof of Theorems 2.3 and 2.4. It is shown in Proposition 4.1 that if $v_P > 0$, the process $x_n = (\theta^{\xi_n}\omega, Y_n)$ is stationary and ergodic under the law induced by the initial distribution $Q$ of $x_0$ [introduced in (12)] and the quenched measure $P_\omega$. The converse part of Theorem 2.3 is given by Proposition 4.8. Proposition 4.5 defines a renewal structure for the paths of the random walk $X_n$ which is the key element in the proof of Theorem 2.4. The proof that in an i.i.d. environment the law of $x_n$ converges toward a limiting distribution (Theorem 2.4) is completed at the end of the section.

Let $\overline{\omega}_n = \theta^{\xi_n}\omega$ and consider the process $x_n = (\overline{\omega}_n, Y_n)$, $n \geq 0$, defined in $(\Omega \times \mathcal{D}, \mathcal{F} \otimes \mathcal{S}_\mathcal{D})$, where $\mathcal{S}_\mathcal{D}$ is the set of all subsets of $\mathcal{D}$. For any $\mu \in \mathcal{M}_d$, $(x_n)_{n\geq 0}$ is a Markov chain under $\mathbb{P}^\mu$ with transition kernel

$$K(\omega, i; B, j) = p_0(i,j)\mathbf{I}_B(\theta\omega) + r_0(i,j)\mathbf{I}_B(\omega) + q_0(i,j)\mathbf{I}_B(\theta^{-1}\omega).$$

By Lemmas 3.1 and 3.2, $\mathbb{E}^\pi(T_1 < \infty) = 1$ if $v_P > 0$. Hence $Q$ is a probability measure. It is shown in Proposition 4.1 below that $Q$ is an invariant distribution for the kernel $K$, that is, $\int_{\Omega \times D} Q(dx)K(x, A) = Q(A)$ for $A \in \mathcal{F} \otimes \mathcal{S}_\mathcal{D}$.

Define a probability measure $\overline{Q}$ on $(\Omega, \mathcal{F})$ by setting

$$\overline{Q}(B) = Q(B, \mathcal{D}), \qquad B \in \mathcal{F},$$

and define a collection $\nu = (\nu_\omega)_{\omega \in \Omega} \in \mathcal{M}_d$ of probability measures on $\mathcal{D}$ by setting

$$\nu_\omega(i) = \frac{dQ_i(\omega)}{d\overline{Q}(\omega)}, \qquad \text{where } Q_i(B) := Q(B, i) \qquad \text{for } B \in \mathcal{F}.$$

For $m \leq 0$ and $i \in \mathcal{D}$ define

$$(26) \qquad N_m^i := \{\#n \in [0, T_1) : \xi_n = m, Y_n = i\}.$$

Note that for any bounded measurable function $f : \Omega \to \mathbb{R}$ and all $i \in \mathcal{D}$,

$$\int_\Omega f(\omega)Q(d\omega, i) = v_P \sum_{n=0}^{\infty} \mathbb{E}^\pi(f(\overline{\omega}_n); Y_n = i, T_1 > n)$$

$$
\begin{aligned}
(27) \qquad &= v_{_P} \sum_{m \le 0} \mathbb{E}^{\pi}(f(\theta^m \omega) N_m^i) \\
&= v_{_P} E_P\left( \sum_{k \in \mathcal{D}} \pi_\omega(k) \sum_{m \le 0} f(\theta^m \omega) E_\omega^k(N_m^i) \right) \\
&= v_{_P} E_P\left( f(\omega) \sum_{k \in \mathcal{D}} \sum_{m \le 0} \pi_{\theta^{-m}\omega}(k) E_{\theta^{-m}\omega}^k(N_m^i) \right).
\end{aligned}
$$

Therefore,

$$
\nu_\omega(i) = \frac{\sum_{k \in \mathcal{D}} \sum_{m \le 0} \pi_{\theta^{-m}\omega}(k) E_{\theta^{-m}\omega}^k(N_m^i)}{\sum_{k \in \mathcal{D}} \sum_{m \le 0} \pi_{\theta^{-m}\omega}(k) \sum_{j \in \mathcal{D}} E_{\theta^{-m}\omega}^k(N_m^j)}, \qquad P\text{-a.s.}
$$

Note that

$$
(28) \qquad\qquad P_\omega(\nu_\omega(i) > 0) = 1 \qquad \forall i \in \mathcal{D},
$$

because $E_\omega^i(N_0^i) \ge 1$.

Set $\mathbb{Q} = \overline{Q} \otimes P_\omega^\nu$. In words, under $\mathbb{Q}$ the environment $\omega$ is distributed according to the measure $\overline{Q}$ while the initial position $Y_0$ of the random walk in environment $\omega$ is determined according to the measure $\nu_\omega$. With a slight abuse of notation we will denote by $\mathbb{Q}$ also the annealed measure induced on the path space of the random walk, that is,

$$
\mathbb{Q}((X_n)_{n \ge 0} \in \cdot) = E_{\overline{Q}}(P_\omega^\nu((X_n)_{n \ge 0} \in \cdot)),
$$

where $E_{\overline{Q}}$ stands for the expectation with respect to $\overline{Q}$.

The following proposition is standard for one-dimensional processes in random environment (see, e.g., [20], [6], Lecture 1, [27], Section 2.1). Note that we do not assume here yet that the environment is an i.i.d. sequence.

PROPOSITION 4.1. *Let Conditions C and D hold and assume in addition that $v_{_P} > 0$. Then:*

(a) *The measure $Q$ is invariant under the kernel $K$, that is,*

$$
Q(B, i) = \sum_{j \in \mathcal{D}} \int_\Omega \int_\Omega \mathbf{I}_B(\omega') K(\omega, j; d\omega', i) Q(d\omega, j).
$$

*In particular, the sequence $(\overline{\omega}_n, Y_n)_{n \ge 0}$ is stationary under the law $\mathbb{Q}$.*

(b) *$\overline{Q} \sim P$ and $\Lambda(\omega) := \frac{d\overline{Q}}{dP} \le dv_{_P} \| u_0 + \sum_{i=1}^\infty a_i a_{i-1} \cdots a_2 a_1 (\mathbf{1} + u_0) \|$, where the vector $u_0$ is defined in (22).*

(c) *The sequence $(\overline{\omega}_n, Y_n)_{n \ge 0}$ is ergodic under the law $\mathbb{Q}$.*

(d) *$Q$ is the unique invariant probability measure for the kernel $K$ which is absolutely continuous with respect to $P \otimes \{$counting measure on $\mathcal{D}\}$.*

PROOF.    (a) The proof is by a "cyclic argument." We have

$$Q(B, i) = v_P \sum_{n=0}^{\infty} \mathbb{P}^\pi(\overline{\omega}_n \in B, Y_n = i; T_1 > n)$$

and, on the other hand,

$$\frac{1}{v_P} QK(B, i) = \frac{1}{v_P} \sum_{j \in \mathcal{D}} \int_\Omega K(\omega, j; B, i) Q(d\omega, j)$$

$$= \sum_{n=0}^{\infty} \mathbb{E}^\pi(\mathbf{I}_B(\overline{\omega}_{n+1}) \mathbf{I}_{\{i\}}(Y_{n+1}); T_1 > n)$$

$$= \sum_{n=0}^{\infty} \mathbb{E}^\pi(\mathbf{I}_B(\overline{\omega}_{n+1}) \mathbf{I}_{\{i\}}(Y_{n+1}); T_1 = n+1)$$

$$+ \sum_{n=0}^{\infty} \mathbb{E}^\pi(\mathbf{I}_B(\overline{\omega}_{n+1}) \mathbf{I}_{\{i\}}(Y_{n+1}); T_1 > n+1)$$

$$= \mathbb{P}^\pi(\overline{\omega}_{T_1} \in B, Y_{T_1} = i) + \sum_{n=1}^{\infty} \mathbb{P}^\pi(\overline{\omega}_n \in B, Y_n = i; T_1 > n).$$

Using (9) and the translation invariance of $P$ with respect to the shift $\theta$,

$$\mathbb{P}^\pi(\overline{\omega}_{T_1} \in B, Y_{T_1} = i) = E_P\left(\sum_{j \in \mathcal{D}} \pi_\omega(j) P_\omega^j(Y_{T_1} = i) \mathbf{I}_B(\theta\omega)\right)$$

$$= E_P(\pi_{\theta\omega}(i) \mathbf{I}_B(\theta\omega)) = E_P(\pi_\omega(i) \mathbf{I}_B(\omega))$$

$$= \mathbb{P}^\pi(Y_0 = i, \overline{\omega}_0 \in B),$$

completing the proof.

(b) It follows from (27) that

$$\Lambda(\omega) = v_P \sum_{k \in \mathcal{D}} \sum_{m \leq 0} \pi_{\theta^{-m}\omega}(k) E_{\theta^{-m}\omega}^k(N_m \cdot \mathbf{1}),$$

where the vectors $N_m = (N_m^1, \ldots, N_m^d)$ are defined in (26).

For $m \leq 0$ and $n \leq 0$, let

$$(29) \qquad M_{m,n}^k := E_\omega(N_m \cdot \mathbf{1} | X_0 = (n, k)),$$

and define $M_{m,n} = (M_{m,n}^1, \ldots, M_{m,n}^d)$.

Conditioning on the first step of the random walk and setting $M_{m,1}^k := 0$ for any $k \in \mathcal{D}$, we obtain for $m + 1 \leq n \leq 0$,

$$M_{m,n} = q_n \mathbf{1}_{\{n = m+1\}} + p_n M_{m,n+1} + r_n M_{m,n} + q_n M_{m,n-1}.$$

For $m+1 \leq n \leq 1$ let $g_{m,n} = M_{m,n-1} - \eta_{n-1}M_{m,n}$. Similarly to (21), we get

$$g_{m,n+1} = a_n \mathbf{1}_{\{n=m+1\}} + a_n g_{m,n}.$$

Recall the vectors $u_n$ defined in (17). Taking into account that $g_{m,n} \leq M_{m,n-1} \leq u_{n-1}$, we obtain for $m \leq -1$,

$$(30) \qquad M_{m,0} = g_{m,1} = a_0 a_{-1} \cdots a_{m+2}(a_{m+1} + a_{m+1}g_{m,m+1})$$

$$\leq a_0 a_{-1} \cdots a_{m+2} a_{m+1}(\mathbf{1} + u_m).$$

It follows that for $m \leq -1$,

$$E^k_{\theta^{-m}\omega}(N_m \cdot \mathbf{1}) = M^k_{m,0}(\theta^{-m}\omega) \leq a_{-m}a_{-m-1}\cdots a_2 a_1(\mathbf{1} + u_0)(k),$$

and hence

$$\Lambda(\omega) \leq v_P \sum_{k \in \mathcal{D}} \left[ \pi_\omega(k)u_0(k) + \sum_{i=1}^{\infty} \pi_{\theta^i\omega}(k)a_i a_{i-1}\cdots a_2 a_1(\mathbf{1} + u_0)(k) \right]$$

$$\leq dv_P \left\| u_0 + \sum_{i=1}^{\infty} a_i a_{i-1} \cdots a_2 a_1(\mathbf{1} + u_0) \right\|.$$

This completes the proof of part (b) of the proposition.

REMARK 4.2. The estimate for $\Lambda(\omega)$ that we obtained is not tight because $u_m$ is a rather rough upper bound for $g_{m,m+1}$ in (30). This estimate, or more precisely, the upper bound (30) for $M_{m,0}$ will be used in the proof of Lemma 5.4.

(c) Our aim is to show that the sequence $(x_n)_{n \geq 0}$ is ergodic under $\mathbb{Q}$. The proof is standard (see, e.g., [6], pages 9–15, or [27], Section 2) and is based on the application of Lemma 3.5. We outline it here because the final step differs slightly from the corresponding argument for RWRE on $\mathbb{Z}^d$.

Let $\mathcal{X} := (x_n)_{n \geq 0}$ and let $A \in \sigma(x_n : n \geq 0)$ be an invariant set, that is, $\eta^{-1}A = A$, $\mathbb{Q}$-a.s., where $\eta$ is the usual shift on the space $(\Omega \times \mathcal{D})^{\mathbb{Z}_+}: (\eta \mathcal{X})_n = x_{n+1}$. It suffices to show that $\mathbb{Q}(A) \in \{0,1\}$. For $x \in \Omega \times \mathcal{D}$ let $\mathbb{P}_x$ be the law of the Markov chain $(x_n)_{n \in \mathbb{Z}_+}$ with the initial state $x_0 = x$. Set $h(x) = \mathbb{P}_x(A)$.

By Lemma 3.5 there exists a set $B \in \mathcal{F} \otimes \mathcal{S}_{\mathcal{D}}$ such that $h(x) = \mathbf{I}_B(x)$, $\mathbb{Q}$-a.s. It remains to show that $Q(B) \in \{0,1\}$. By the martingale property of the sequence $h(x_n)$, $\mathbf{I}_B(x) = K\mathbf{I}_B(x)$, $\mathbb{Q}$-a.s. Write $B = \bigcup_{i \in \mathcal{D}} B_i \times \{i\}$. We have $P$-a.s. (since $\overline{Q} \sim P$) for all $i \in \mathcal{D}$,

$$\mathbf{I}_{B_i}(\omega) = \mathbf{I}_B(\omega, i) = K\mathbf{I}_B(\omega, i)$$

$$= \sum_{j \in \mathcal{D}} [p_0(i,j)\mathbf{I}_{B_j}(\theta\omega) + r_0(i,j)\mathbf{I}_{B_j}(\omega) + q_0(i,j)\mathbf{I}_{B_j}(\theta^{-1}\omega)].$$

Since $\mathbf{I}_{B_i}$ takes only two values, either 0 or 1, $\sum_{j \in \mathcal{D}} p(i,j) + r(i,j) + q(i,j) = 1$, and $p_0(i,j) > 0$ for any $i, j \in \mathcal{D}$, $P$-a.s., we obtain that

$$(31) \qquad \mathbf{I}_{B_i}(\omega) = \mathbf{I}_{B_j}(\theta\omega) \qquad \text{for all } i, j \in \mathcal{D}, \qquad P\text{-a.s.}$$

In particular, $\mathbf{I}_{B_i}(\omega)$ are invariant functions of the ergodic sequence $(\omega_n)_{n \in \mathbb{Z}}$ and hence $P(B_i) \in \{0, 1\}$. Furthermore, by (31) and the translation invariance of the measure $P$, $P(B_i) = P(\theta B_j) = P(B_j)$ for every $i, j \in \mathcal{D}$. It follows from (31) that

$$P \otimes \mathrm{lm}_d(B) \in \{0, 1\},$$

where $\mathrm{lm}_d$ is the counting probability measure on $\mathcal{D}$, that is, $\mathrm{lm}_d(i) = 1/d$ for all $i \in \mathcal{D}$. This completes the proof of claim (c) because $\overline{Q} \sim P$ and hence, by (28), $Q \sim P \otimes \mathrm{lm}_d$.

REMARK 4.3.   Note that the only property of the stationary measure $\mathbb{Q}$ that we use in the proof is that $\mathbb{Q} = \overline{Q} \otimes P_\omega^\nu$ and $\overline{Q} \otimes \nu_\omega \sim P \otimes \mathrm{lm}_d$. This observation will be used in the proof of Proposition 4.8.

(d) Let $Q_0$ be an invariant probability measure for the kernel $K$ such that $Q_0 \ll P \otimes \mathrm{lm}_d$. It follows from the ergodic theorem that for any measurable set $A \in \mathcal{F} \otimes \mathcal{S}_\mathcal{D}$, we have $\overline{Q} \otimes \pi_\omega$-a.s.,

$$(32) \qquad \frac{1}{n}\sum_{i=0}^{n-1}\mathbb{Q}_x(x_i \in A) = \frac{1}{n}\sum_{i=0}^{n-1} K^i(x, A) \underset{n \to \infty}{\to} \mathbb{Q}(A).$$

Taking the integral with respect to $Q_0(dx)$, and using the fact that $\overline{Q} \sim P$, we obtain that $Q_0(A) = \mathbb{Q}(A)$. The proof of the proposition is completed. $\square$

We turn now to the proof that *in an i.i.d. environment* $x_n$ converges under $\mathbb{P}^\mu$ to the limit distribution $Q$. The proof is similar to that of the analogous statement for RWRE on $\mathbb{Z}^d$ (see [17] for $d = 1$ and [26] for $d > 1$) and uses a renewal structure that we proceed now to introduce.

For $n \in \mathbb{N}$, let

$$(33) \qquad \chi_n = \inf_{i > n} \xi_i,$$

and fix any $i^* \in \mathcal{D}$ such that

$$(34) \qquad \mathbb{P}^{i^*}(\chi_0 \leq 0) < 1.$$

Such $i^*$ exists because $\mathbb{P}^i(\chi_0 \leq 0) = 1$ for all $i \in \mathcal{D}$ means (by the strong Markov property of the random walk in a fixed environment) that $\mathbb{P}^\mu(\xi_n \leq 0 \text{ i.o.}) = 1$ for all $\mu \in \mathcal{M}_d$. This would imply that $X_n$ is either transient to the left or recurrent and hence contradicts Condition D.

Let $T_0^* = 0$ and for $n \geq 1$ define

$$T_n^* = \inf\{m > T_{n-1}^* : m = T_k \text{ for some } k \geq 0 \text{ and } Y_m = i^*\}.$$

Thus $T_n^*$ is a subsequence of random times $T_k$ corresponding to the indexes $k$ such that $Y_{T_k} = i^*$, that is, the random walk hits the layer $\{k\} \times \mathcal{D}$ at height $i^*$.

Next, let $\rho_{-1} = -1$ and for $n \geq 0$ define

$$\rho_n = \inf\left\{m > \rho_{n-1} : \max_{j \leq m} \xi_j = \xi_m < \min_{j > m} \xi_j \text{ and } Y_m = i^*\right\}.$$

Note that $(\rho_n)_{n \geq 0}$ is a subsequence of $(T_n^*)_{n \geq 0}$. At times $\rho_n$: (i) the path of the random walk hits some layer $\{m\} \times \mathcal{D}$ at first time, (ii) its altitude (the second component) is $i^*$, and (iii) the random walk will never return again to this layer.

LEMMA 4.4. *Let the conditions of Theorem 2.4 hold. Then for any $\mu \in \mathcal{M}_d$ we have $\mathbb{P}^\mu(\rho_n < \infty \text{ for all } n \geq 0) = 1$.*

PROOF. First, we shall prove that

(35) $$\mathbb{P}^\mu(T_n^* < \infty \text{ for all } n \in \mathbb{N}) = 1.$$

Toward this end observe that for any $k \in \mathbb{N}$ and any increasing sequence of positive integers $(\alpha_n)_{n \in \mathbb{N}}$ we have

$$\begin{aligned}
\mathbb{P}^\mu(T_k^* < \infty, T_{k+1}^* = \infty) &\leq \sum_{n \in \mathbb{N}} E_P(P_\omega^\mu(T_k^* = T_n, T_{k+1}^* > \alpha_n)) \\
&= \sum_{n \in \mathbb{N}} E_P(P_\omega^\mu(T_k^* = T_n) P_{\theta^n \omega}^{i^*}(T_1^* > \alpha_n)) \\
&\leq \sum_{n \in \mathbb{N}} \mathbb{P}^{i^*}(T_1^* > \alpha_n).
\end{aligned}$$

(36)

With $\eta_n^*(i,j) := \eta_n(i,j) \mathbf{I}_{\{j \neq i^*\}}$, we have for any $n \in \mathbb{N}$,

$$P_\omega^\mu(T_1^* = \infty) \leq P_\omega^\mu(T_1^* > T_n) \leq \|\eta_0^* \cdots \eta_{n-1}^*\| \leq \|\eta_0^*\| \cdots \|\eta_{n-1}^*\|.$$

By condition C4 and the ergodic theorem,

$$\lim_{n \to \infty} \frac{1}{n} \sum_{k=0}^{n-1} \log \|\eta_k^*\| = E_P(\log \|\eta_0^*\|) < 0, \qquad P\text{-a.s.}$$

Therefore, $\mathbb{P}^\mu(T_1^* = \infty) = \lim_{n \to \infty} \mathbb{P}^\mu(T_1^* > T_n) = 0$ by the bounded convergence theorem. This shows that we can choose a sequence $\alpha_n$ to make the rightmost expression in (36) arbitrarily small and therefore proves (35).

Let $A$ denote the event $\{\rho_n = \infty$ for some $n \geq 0\}$. Our aim is to prove that $\mathbb{P}^\mu(A) = 0$ for any $\mu \in \mathcal{M}_d$. For $n \geq 0$, let

$$k_n = \inf\{i > T_n^* : \xi_i = \xi_{T_n^*}\}.$$

Then

$$A = \bigcup_{m \geq 0} A_m \qquad \text{where } A_m := \{k_n < \infty \text{ for all } n \geq m\}.$$

Fix now any $m \geq 0$, let $h_0 = \xi_{T_m^*}$, $g_0 = m$, and on the event $A_m$ define for $n \geq 0$,

$$h_{n+1} = \inf\{j > k_{g_n} : j = T_i^* \text{ for some } i \in \mathbb{N}\} = T_{g_{n+1}}^*,$$

where the last equality defines random variables $g_n$, $n \geq 1$. We note that this construction is reminiscent of the renewal structure for RWRE on $\mathbb{R}^d$ introduced in [26] and $h_n$ are a counterpart of their "fresh times." In words, $h_{n+1} = T_{g_{n+1}}^*$ is the first time after $k_{g_n}$ (which is the time of the first return to level $\xi_{T_{g_n}^*}$ after $T_{g_n}^*$) when the path of the random walk hits a fresh part of the environment (i.e., a new level on the right of $\xi_{T_{g_n}^*}$) and its height at that moment is $i^*$.

For the fixed value of the parameter $m$ and $n \in \mathbb{N}$, let $A_{m,n}$ denote the event $\{k_{g_i} < \infty$ for $i = 0, \ldots, n\}$. Then, since $P$ is a product measure,

$$\mathbb{P}^\mu(A_m) \leq \mathbb{P}^\mu(A_{m,n+1}) = \sum_{x \in \mathbb{Z}} \mathbb{P}^\mu\left(A_{m,n}, \xi_{h_n} = x, \inf_{j > h_n} \xi_j \leq x\right)$$

$$= \sum_{x \in \mathbb{Z}} E_P(P_\omega^\mu(A_{m,n}, \xi_{h_n} = x) P_{\theta^x \omega}^{i^*}(\chi_0 \leq 0))$$

$$= \sum_{x \in \mathbb{Z}} \mathbb{P}^\mu(A_{m,n}, \xi_{h_n} = x) \mathbb{P}^{i^*}(\chi_0 \leq 0)$$

$$= \mathbb{P}^\mu(A_{m,n}) \mathbb{P}^{i^*}(\chi_0 \leq 0) \leq \cdots \leq [\mathbb{P}^{i^*}(\chi_0 \leq 0)]^n \underset{n \to \infty}{\to} 0,$$

where we use the fact that the random variable $P_\omega^\mu(A_{m,n}, \xi_{h_n} = x)$ is measurable with respect to $\sigma(\omega_n : n < x)$ while $P_{\theta^x \omega}^{i^*}(\chi_0 \leq 0)$ is measurable with respect to $\sigma(\omega_n : n \geq x)$. The proof of the lemma is therefore complete. $\quad\square$

Define a probability measure $\widetilde{\mathbb{P}}$ on $(\mathbb{Z} \times \mathcal{D})^{\mathbb{Z}_+} \times \Omega$ by setting

$$(37) \qquad \widetilde{\mathbb{P}}((X_n)_{n \geq 0} \in \cdot, \omega \in \cdot) = \mathbb{P}^{i^*}((X_n)_{n \geq 0} \in \cdot, \omega \in \cdot | \chi_0 > 0).$$

Note that for any $A \in \sigma(X_n : n \geq 0)$,

$$\widetilde{\mathbb{P}}(A) = \frac{\mathbb{P}^{i^*}(A \cap \chi_0 > 0)}{\mathbb{P}^{i^*}(\chi_0 > 0)} = \frac{E_P(P_\omega^{i^*}(A \cap \chi_0 > 0))}{E_P(P_\omega^{i^*}(\chi_0 > 0))}$$

$$(38) \qquad = \frac{E_P(\pi_\omega(i^*)) \cdot E_P(P_\omega^{i^*}(A \cap \chi_0 > 0))}{E_P(\pi_\omega(i^*)) \cdot E_P(P_\omega^{i^*}(\chi_0 > 0))}$$

$$= \frac{E_P(\pi_\omega(i^*) P_\omega^{i^*}(A \cap \chi_0 > 0))}{E_P(\pi_\omega(i^*) P_\omega^{i^*}(\chi_0 > 0))}$$

$$= \frac{\mathbb{P}^\pi(A \cap Y_{T_0} = i^* \cap \chi_0 > 0)}{\mathbb{P}^\pi(Y_{T_0} = i^* \cap \chi_0 > 0)} \mathbb{P}^\pi(A | Y_{T_0} = i^*, \chi_0 > 0),$$

where we use the fact that the random vector $\pi_\omega$ is measurable with respect to $\sigma(\omega_n : n < 0)$ while $P_\omega^{i^*}(A \cap \chi_0 > 0)$ is measurable with respect to $\sigma(\omega_n : n \geq 0)$.

For $k \geq 0$, denote

$$\mathcal{G}_k = \sigma(\rho_1, \dots, \rho_k, X_0, \dots, X_{\rho_k}, (\omega_n)_{n < \xi_{\rho_k}}).$$

The next proposition is the key to the proof of Theorem 2.4.

PROPOSITION 4.5. *Let the conditions of Theorem 2.4 hold. Then,*

$$(39) \qquad \mathbb{P}^\mu((X_{\rho_k+n} - X_{\rho_k})_{n \geq 0} \in \cdot, (\omega_{\xi_{\rho_k}+i})_{i \geq 0} \in \cdot | \mathcal{G}_k)$$

$$= \widetilde{\mathbb{P}}((X_n)_{n \geq 0} \in \cdot, (\omega_i)_{i \geq 0} \in \cdot)$$

*for any $\mu \in \mathcal{M}_d$ and $k \geq 0$.*

COROLLARY 4.6. *Under the conditions of Theorem 2.4, random $\mathbb{R}^2$-vectors $(\xi_{\rho_0}, \rho_0)$, $(\xi_{\rho_1} - \xi_{\rho_0}, \rho_1 - \rho_0), \dots, (\xi_{\rho_n} - \xi_{\rho_{n-1}}, \rho_n - \rho_{n-1}), \dots$ are independent under $\mathbb{P}^\mu$ for any $\mu \in \mathcal{M}_d$. Furthermore, the random vectors $(\xi_{\rho_1} - \xi_{\rho_0}, \rho_1 - \rho_0), \dots, (\xi_{\rho_n} - \xi_{\rho_{n-1}}, \rho_n - \rho_{n-1}), \dots$ are distributed under $\mathbb{P}^\mu$ as $(\xi_{\rho_1}, \rho_1)$ under $\widetilde{\mathbb{P}}$.*

The proof of the proposition and the corollary is similar to the proof of Theorem 1.4 in [26] (see also Lemma 3.2.3 in [27]). For the sake of completeness we outline it here, following the proof in [27], pages 263–264. First, observe that for any integer number $m \geq k$ and $\mathcal{G}_k$-measurable random variable $h$ there exists a random variable $h_m$, measurable with respect to $\sigma((X_i - X_0)_{i \leq T_m}, (\omega_n)_{n < m})$, such that $h \mathbf{I}_{\{\rho_k = T_m\}} = h_m \mathbf{I}_{\{\chi_m > m\}}$. Therefore, for any measurable sets $A \subset (\mathbb{R}^2)^{\mathbb{Z}_+}$ and $B \subset \mathcal{S}^{\mathbb{Z}_+}$ we have

$$\mathbb{E}^\mu(\mathbf{I}_A[(X_{\rho_k+n} - X_{\rho_k})_{n \geq 0}] \cdot \mathbf{I}_B[(\omega_{\xi_{\rho_k}+i})_{i \geq 0}] \cdot h)$$

$$= \sum_{m \geq k} E_P(\mathbf{I}_B[(\omega_{m+i})_{i \geq 0}] \cdot E_\omega^\mu(\mathbf{I}_A[(X_{T_m+n} - X_{T_m})_{n \geq 0}] \cdot h \cdot \mathbf{I}_{\{\rho_k = T_m\}}))$$

$$= \sum_{m \geq k} E_P(\mathbf{I}_B[(\omega_{m+i})_{i \geq 0}] \cdot E_{\theta^m \omega}^{i^*}(\mathbf{I}_A[(X_n - X_0)_{n \geq 0}] \cdot \mathbf{I}_{\{\chi_0 > 0\}}) \cdot E_\omega^\mu(h_m))$$

$$= \sum_{m \geq k} \mathbb{E}^{i^*}(\mathbf{I}_A[(X_n - X_0)_{n \geq 0}] \cdot \mathbf{I}_B[(\omega_i)_{i \geq 0}] \cdot \mathbf{I}_{\{\chi_0 > 0\}}) \cdot \mathbb{E}^\mu(h_m),$$

where we use the Markov property in the second equality and the i.i.d. structure of the environment in the last one. Substituting $A = (\mathbb{R}^2)^{\mathbb{Z}_+}$ and $B = \mathcal{S}^{\mathbb{Z}_+}$ yields

$$\mathbb{E}^\mu(h) = \sum_{m \geq k} \mathbb{E}^{i^*}(\mathbf{I}_{\{\chi_0 > 0\}}) \cdot \mathbb{E}^\mu(h_m) = \mathbb{P}^{i^*}(\chi_0 > 0) \cdot \sum_{m \geq k} \mathbb{E}^\mu(h_m).$$

Therefore, we get

$$\mathbb{E}^\mu(\mathbf{I}_A[(X_{\rho_k + n} - X_{\rho_k})_{n \geq 0}] \cdot \mathbf{I}_B[(\omega_{\xi_{\rho_k} + i})_{i \geq 0}] \cdot h)$$
$$= \widetilde{\mathbb{E}}(\mathbf{I}_A[(X_n - X_0)_{n \geq 0}] \cdot \mathbf{I}_B[(\omega_i)_{i \geq 0}]) \cdot \mathbb{E}^\mu(h),$$

which is equivalent to (39). The corollary follows from (39) by considering the finite-dimensional distributions and using induction.

The next step toward the proof of Theorem 2.4 is to show that $\widetilde{\mathbb{E}}(\rho_1) < \infty$, where $\widetilde{\mathbb{E}}$ is the expectation with respect to the measure $\widetilde{\mathbb{P}}$ defined in (37).

LEMMA 4.7. *Let the conditions of Theorem 2.4 hold. Then $\widetilde{\mathbb{E}}(\rho_1) < \infty$.*

PROOF. Recall the random variables $\chi_n$ defined in (33) and let random numbers $\lambda_n$ be defined by the equation $\rho_n = T_{\lambda_n}$, $n \geq 0$. The sequence $(Y_{T_n}, \chi_n)_{n \in \mathbb{N}}$ is stationary under $\mathbb{P}^\pi$, and the definition of $\lambda_n$ can be rewritten on the event $\{Y_0 = i^*, \chi_0 > 0\}$ as $\lambda_n = \inf\{m > \lambda_{n-1} : Y_{T_m} = i^*, \chi_m > m\}$ with $\lambda_0 = 0$. That is, $\lambda_n$ for $n > 0$ can be viewed as the time of $n$th return to the set $\{i^*\} \times \mathbb{N}$ by the sequence $(Y_{T_m}, \chi_m - m)_{m \geq 0}$. By Lemma 4.4, $\mathbb{P}^\pi(Y_{T_k} = i^*$ and $\chi_k > 0$ i.o.$) = 1$. Therefore, by Kac's recurrence lemma (cf. Theorem 3.3 in [12], page 348),

$$\widetilde{\mathbb{E}}(\lambda_1) = \mathbb{E}^\pi(\lambda_1 | Y_{T_0} = i^*, \chi_0 > 0) = 1/\mathbb{P}^\pi(Y_{T_0} = i^*, \chi_0 > 0) < \infty,$$

where the first equality is due to the identity (38).

On one hand, it follows from Corollary 4.6 that

$$\rho_n/n = T_{\lambda_n}/n \underset{n \to \infty}{\to} \widetilde{\mathbb{E}}(\rho_1), \qquad \widetilde{\mathbb{P}}\text{-a.s.},$$

and on the other hand, by Lemma 3.2 and Corollary 4.6,

$$\frac{T_{\lambda_n}}{n} = \frac{T_{\lambda_n}}{\lambda_n} \cdot \frac{\lambda_n}{n} \underset{n \to \infty}{\to} \mathbb{E}^\pi(\tau_1) \cdot \widetilde{\mathbb{E}}(\lambda_1), \qquad \widetilde{\mathbb{P}}\text{-a.s.}$$

Therefore, $\widetilde{\mathbb{E}}(\rho_1) = \mathbb{E}^\pi(\tau_1) \cdot \widetilde{\mathbb{E}}(\lambda_1) < \infty$ and the proof of the lemma is completed. □

We are now in the position to complete the proof of Theorem 2.4. The rest of the proof is similar to that given in [17] and [26] and relies on an application of the renewal theorem.

COMPLETION OF THE PROOF OF THEOREM 2.4. By the ergodic theorem and Proposition 4.1, for any bounded measurable function $f : \Omega \times \mathcal{S}_{\mathcal{D}} \to \mathbb{R}$ and the Markov chain $x_n = (\theta^{\xi_n} \omega, Y_n)$,

$$\lim_{n \to \infty} \frac{1}{n} \sum_{i=0}^{n-1} f(x_i) = \mathbb{E}_{\mathbb{Q}}(f(x_0)), \qquad \mathbb{Q}\text{-a.s.},$$

where $\mathbb{E}_{\mathbb{Q}}$ stands for the expectation with respect to the measure $\mathbb{Q}$. By part (ii) of Proposition 4.1 and (28), the convergence is also $\mathbb{P}^{\mu}$-a.s. Hence, if $\lim_{n \to \infty} \mathbb{E}^{\mu}(f(x_n))$ exists it is equal to $\mathbb{E}_{\mathbb{Q}}(f(x_0))$. Therefore, it suffices to prove the existence of the limit for arbitrary bounded continuous $f$.

Let $f : \Omega \times \mathcal{D} \to \mathbb{R}$ be a bounded continuous function which depends on finitely many coordinates only, that is, $f$ is measurable with respect to $\sigma(\omega_n : |n| < N) \otimes \mathcal{S}_{\mathcal{D}}$ for some $N \in \mathbb{N}$. For any $i \in \mathcal{D}$ and $n \geq 1$, write

$$\mathbb{E}^{\mu}(f(\overline{\omega}_n, Y_n)) = \mathbb{E}^{\mu}(f(\overline{\omega}_n, Y_n); \rho_{N+1} > n) + \mathbb{E}^{\mu}(f(\overline{\omega}_n, Y_n); \rho_{N+1} \leq n).$$

The first term in the right-hand side tends to zero when $n$ goes to infinity. To evaluate the second term we write

$$\mathbb{E}^{\mu}(f(\overline{\omega}_n, Y_n; \rho_{N+1} \leq n))$$

$$= \sum_{i=1}^{\infty} \mathbb{E}^{\mu}(f(\overline{\omega}_n, Y_n); \rho_{i+N} \leq n < \rho_{i+N+1})$$

$$= \sum \mathbb{E}^{\mu}(f(\theta^y \omega, j); \rho_{i+N} \leq n < \rho_{i+N+1}, \rho_i = m, \xi_{\rho_i} = x, X_n = (y, j)),$$

where the last sum is taken over all $i, m \in \mathbb{N}, x, y \in \mathbb{Z}, j \in \mathcal{D}$. Note that on the event $\{\rho_{i+N} \leq n < \rho_{i+N+1}, \rho_i = m, \xi_{\rho_i} = x, X_n = (y, j)\}$, the random variable $f(\theta^y \omega, j)$ is $\sigma(\omega_z : z \geq x)$-measurable because, by the definition of the random times $\rho_i$,

$$x = \xi_{\rho_i} \leq \xi_{\rho_{i+N}} - N \leq y - N.$$

Therefore, by Proposition 4.5,

$$\mathbb{E}^{\mu}(f(\overline{\omega}_n, Y_n); \rho_{N+1} \leq n)$$

$$= \sum_{i=1}^{\infty} \sum_{m \in \mathbb{N}} \sum_{x, y \in \mathbb{Z}} \sum_{j \in \mathcal{D}} \mathbb{P}^{\mu}(\rho_i = m, \xi_{\rho_i} = x)$$

$$\times \widetilde{\mathbb{E}}(f(\theta^{y-x} \omega, j); \rho_N \leq n - m < \rho_{N+1}, X_{n-m} = (y - x, j))$$

$$= \sum_{t=0}^{n} \mathbb{P}^{\mu}(n - t = \rho_i \text{ for some } i \in \mathbb{N}) \cdot \widetilde{\mathbb{E}}(f(\overline{\omega}_t, Y_t); \rho_N \leq t < \rho_{N+1}).$$

By condition C4 and (34), $\mathbb{P}^{i^*}(\rho_1 = 1) > 0$. It follows from Corollary 4.6 and the renewal theorem for arithmetic distributions (cf. [13], Volume II, page

347) that

$$\mathbb{P}^{\mu}(n - t = \rho_i \text{ for some } i \in \mathbb{N})$$

$$= \mathbb{P}^{\mu}(n - t - \rho_0 = \rho_i - \rho_0 \text{ for some } i \in \mathbb{N}) \underset{n \to \infty}{\to} \frac{1}{\widetilde{\mathbb{E}}(\rho_1)}.$$

Note that by Proposition 4.5, random variable $n - t - \rho_0$ is independent of the sequence $(\rho_n - \rho_{n-1})_{n \geq 1}$. Therefore, one can first condition on the value of "delay" $\rho_0$, apply the cited renewal theorem, and then use the bounded convergence theorem in order to get this result.

By Proposition 4.5 and Lemma 4.7, $\sum_{t=0}^{\infty} \widetilde{\mathbb{P}}(\rho_N \leq t < \rho_{N+1}) = \widetilde{\mathbb{E}}(\rho_1)$ is finite. Therefore, the following limit exists and is a finite number:

$$\lim_{n \to \infty} \mathbb{E}^{\mu}(f(\overline{\varpi}_n, Y_n); \rho_{N+1} \leq n)$$

$$= \frac{1}{\widetilde{\mathbb{E}}(\rho_1)} \sum_{t=0}^{\infty} \widetilde{\mathbb{E}}(f(\overline{\varpi}_t, Y_t); \rho_N \leq t < \rho_{N+1}).$$

This proves the convergence of $\mathbb{E}^{\mu}(f(\overline{\varpi}_n, Y_n))$ for bounded functions $f$ which depend on finitely many coordinates only. Since we can assume without loss of generality that the space $\Omega \times \mathcal{D}$ is compact [because we can consider the sequence $\widetilde{\omega}_n = (p_n, r_n, q_n)$, $n \in \mathbb{Z}$, as an environment], the claim for arbitrary bounded continuous functions follows by a standard argument and we omit it. This completes the proof of the theorem. $\square$

We conclude this section with the proof of the converse part of Theorem 2.3. The method of the proof of the LLN for $\xi_n$ with the help of the auxiliary Markov chain $\hat{x}_n = (x_n, \xi_{n+1} - \xi_n)$ goes back to Kozlov [20] while the idea to use [3], Theorem 1, for the proof that the asymptotic speed of a RWRE is strictly positive is due to Brémont (cf. [7], Theorem 3.5(ii)).

PROPOSITION 4.8.    *Let Conditions C and D hold and suppose in addition that there exists a stationary distribution $Q_0$ for the Markov chain $x_n = (\theta^{\xi_n}\omega, Y_n), n \geq 0$, which is equivalent to $P \otimes \{counting\ measure\ on\ \mathcal{D}\}$. Then $v_P > 0$ and $Q_0$ coincides with the measure $Q$ defined in (12).*

PROOF.    Let $\mathbb{Q}_0$ denote the joint law of the environment $\omega$ and of the random walk $(X_n)_{n \geq 0}$ with $\xi_0 = 0$, induced by the initial distribution $Q_0$ of $(\omega, Y_0)$ and by the quenched measure $P_\omega$ on the path space of the random walk $X_n$. By Remark 4.3, the Markov chain $x_n = (\theta^{\xi_n}\omega, Y_n)$, $n \geq 0$, is stationary and ergodic under $\mathbb{Q}_0$. Consider now the Markov chain $\hat{x}_n = (x_n, \xi_{n+1} - \xi_n)$, $n \geq 0$. Using again Lemma 3.5 and the fact that the transitions of $\hat{x}_n$ depend on the position of $x_n$ only but not on the value of $\xi_n - \xi_{n-1}$, it is not hard to see that the process $(\hat{x}_n)_{n \geq 0}$ is stationary and ergodic under $\mathbb{Q}_0$. Therefore, by the ergodic theorem,

$$\frac{\xi_n}{n} = \frac{1}{n} \sum_{n=1}^{n} (\xi_n - \xi_{n-1}) \underset{n \to \infty}{\to} \mathbb{E}_{\mathbb{Q}_0}(\xi_1), \qquad \mathbb{Q}_0\text{-a.s. and } \mathbb{P}^\pi\text{-a.s.},$$

where $\mathbb{E}_{\mathbb{Q}_0}$ stands for the expectation with respect to the measure $\mathbb{Q}_0$. Moreover, since $X_n$ is transient to the right, $\mathbb{E}_{\mathbb{Q}_0}(\xi_1) \geq 0$. But if $\mathbb{E}_{\mathbb{Q}_0}(\xi_1) = 0$, then $\mathbb{Q}_0(\xi_n = 0 \text{ i.o.}) = 1$ by the recurrence theorem for partial sums of a transformation of stationary and ergodic sequence due to [3]. Since $\mathbb{P}^\pi$ is equivalent to $\mathbb{Q}_0$, this would imply that $\mathbb{P}^\pi(\xi_n = 0 \text{ i.o.}) = 1$ as well, and hence contradicts Condition D. Therefore $v_P > 0$ and $Q_0$ coincides with $Q$ by Proposition 3.2. $\square$

## 5. Annealed CLT for uniformly mixing environments.

This section is devoted to the proof of Theorem 2.5. The annealed CLT for $\xi_n$ is derived from the corresponding result for the hitting times $T_n$ by using the following observation.

LEMMA 5.1. *Let the conditions of Theorem 2.5 hold and assume that for some $\mu \in \mathcal{M}_d$ and suitable constants $\alpha > 0$ and $\sigma > 0$, $(T_n - \alpha n)/\sqrt{n} \Rightarrow N(0, \sigma^2)$ under $\mathbb{P}^\mu$. Then under $\mathbb{P}^\mu$, $(\xi_n - \alpha^{-1}n)/\sqrt{n} \Rightarrow N(0, \sigma^2 \alpha^{-3})$.*

PROOF. The proof is similar to that for RWRE on $\mathbb{Z}$ (see [18] or [2]) except one point where we need to take into account that the sequence $Y_{T_n}$ (i.e., the position of the random walk within the layer $n$ when it visits the layer first time) is not necessarily stationary under $\mathbb{P}^\mu$. For any positive integers $a, b, n$ we have

$$(40) \qquad \{T_a \geq n\} \subset \{\xi_n \leq a\} \subset \{T_{a+b} \geq n\} \cup \left\{\inf_{i \geq T_{a+b}} \xi_i - (a+b) \leq -b\right\}.$$

The first inclusion above expresses the fact that $\xi_n \leq a$ unless the level $a$ has been already crossed by the random walk before time $n$, whereas the second one tells that if $\xi_n \leq a$ but $T_{a+b} < n$, then the maximal excursion of the random walk to the left after $T_{a+b}$ and before escaping to $+\infty$ is at least of the length $(a+b) - a = b$.

Since $X_n$ is transient to the right and

$$\mathbb{P}^\mu\left(\inf_{i \geq T_{a+b}} \xi_i \leq a\right) = E_P\left(\sum_{j \in \mathcal{D}} P^j_{\theta^{a+b}\omega}(\chi_0 \leq -b) P^\mu_\omega(Y_{T_{a+b}} = j)\right)$$

$$\leq \max_{j \in \mathcal{D}} \mathbb{P}^j(\chi_0 \leq -b),$$

the probability of the rightmost event in (40) can be made arbitrarily small uniformly in $n$ and $a$ by fixing $b$ large.

It follows from (40) that

$$(41) \qquad \mathbb{P}^\mu(T_a \geq n) \leq \mathbb{P}^\mu(\xi_n \leq a) \leq \mathbb{P}^\mu(T_{a+b} \geq n) + \max_{j \in \mathcal{D}} \mathbb{P}^j(\chi_0 \leq -b).$$

To complete the proof one can set $a = a(n) = nv_P + xv_P^{3/2}\sqrt{n} + o(\sqrt{n})$, $x \in \mathbb{R}$, and then pass to the limit in (41) using the limit theorem for $T_n$, letting first $n$ and then $b$ go to infinity. The full argument is detailed in [2], page 344, and thus is omitted here. □

In order to show that the CLT holds for $T_n$ we shall apply to $Z_n = \tau_n - v_P^{-1}$ the following general central limit theorem for uniformly mixing sequences (see [12], page 427).

THEOREM 5.2. *Let* $(Z_n)_{n \in \mathbb{Z}}$ *be a stationary random sequence such that* $\mathbb{E}^\pi(Z_0) = 0$ *and* $\mathbb{E}^\pi(Z_0^2) < \infty$. *For* $n \in \mathbb{N}$ *let* $\mathcal{F}^n = \sigma(Z_i : i \geq n)$, $\mathcal{F}_n = \sigma(Z_i : i \leq n)$,

$$(42) \qquad \beta(n) = \sup\{\mathbb{P}^\pi(A|B) - \mathbb{P}^\pi(A) : A \in \mathcal{F}^n, B \in \mathcal{F}_0, \mathbb{P}^\pi(B) > 0\}$$

*and suppose that* $\sum_{n=1}^\infty \sqrt{\beta(2^n)} < \infty$. *If* $S_n = \sum_{i=1}^n Z_i$, *then* $S_n/\sqrt{n} \Rightarrow \sigma N(0,1)$, *where*

$$(43) \qquad \sigma^2 = \mathbb{E}^\pi(Z_1^2) + 2\sum_{n=2}^\infty \mathbb{E}^\pi(Z_1 Z_n).$$

We note that once the convergence of $S_n/\sqrt{n}$ under $\mathbb{P}^\pi$ is established, the convergence under $\mathbb{P}^\mu$ for an arbitrary $\mu \in \mathcal{M}_d$ follows from part (c) of Lemma 3.3 along with assumption (iv) of Theorem 2.5. Indeed, for any $m > 0$, $(T_n - nv_P^{-1})/\sqrt{n} \Rightarrow \sigma N(0,1)$ if and only if $(T_n - T_m - (n-m)v_P^{-1})/\sqrt{n} \Rightarrow \sigma N(0,1)$, and the distribution functions of $(T_n - T_m - (n-m)v_P^{-1})/\sqrt{n}$ under $\mathbb{P}^\pi$ and $\mathbb{P}^\mu$ differ by a factor which is $o(1)$ as $m \to \infty$. Namely, letting $R_{m,n} := T_n - T_m - (n-m)v_P^{-1}$, we have

$$|\mathbb{P}^\mu(R_{m,n}/\sqrt{n} \leq x) - \mathbb{P}^\pi(R_{m,n}/\sqrt{n} \leq x)|$$

$$\leq \sum_{i \in \mathcal{D}} E_P(|P_\omega^\mu(Y_m = i) - P_\omega^\pi(Y_m = i)| \cdot P_\omega^i(R_{0,m-n}/\sqrt{n} \leq x))$$

$$\leq d \cdot E_P(c_0 \cdots c_m P_\omega^i(R_{0,m-n}/\sqrt{n} \leq x)) \leq d \cdot E_P(c_0 \cdots c_m).$$

Taking $m$ to infinity completes the proof of the claim.

Therefore, in order to prove Theorem 2.5 it suffices to check that the second moment condition and the mixing condition of Theorem 5.2 hold for $Z_n = \tau_n - v_P$, and that the limiting variance $\sigma^2$ defined in (43) is strictly positive under the conditions of Theorem 2.5.

The next proposition shows that the sequence $\tau_n - v_P$ is uniformly mixing under $\mathbb{P}^\pi$ with a mixing rate fast enough to apply the general CLT. In the case of one-dimensional RWRE closely related results can be found in [21] (cf. Lemma 4) and [27] (cf. Lemma 2.1.10). The proof for RWRE on strips is more involved but is based on essentially the same idea. Namely, it holds with a large probability that $\inf_{i \geq T_{m+n}} \xi_i \geq m + n/2$, that is, the $\sigma$-fields $\mathcal{F}_m$ and $\mathcal{F}^{m+n}$ are essentially well separated enabling the use of the mixing properties of the environment. Note that by part (iii) of Lemma 3.3, $P_\omega(Y_{T_{m+n}} = j | Y_{T_m} = i) \approx \pi_{\theta^{n+m}\omega}(j)$ uniformly on $i, j \in \mathcal{D}$ and $m \geq 0$ for large $n$, that is, the random walk "forgets" its starting position $Y_0$.

PROPOSITION 5.3. *Let the conditions of Theorem 2.5 hold and let $Z_n = \tau_n - v_P$, $n \in \mathbb{N}$. Then $\sum_{n=1}^\infty \sqrt{\beta(n)} < \infty$ where the mixing coefficients $\beta(n)$ are defined in* (42).

PROOF. For $m \in \mathbb{N}$ and $n \in \mathbb{N} \cup \{\infty\}$ let $\mathcal{T}_m^n$ denote $(Z_i)_{i=m}^n$. For $n, m \in \mathbb{N}$, measurable events $A, B$ and $\mu \in \mathcal{M}_d$, write

$$\mathbb{P}^\mu(\mathcal{T}_{m+3n+1}^\infty \in A | \mathcal{T}_1^m \in B)$$

$$= \frac{1}{\mathbb{P}^\mu(\mathcal{T}_1^m \in B)} \cdot E_P\left(\sum_{i \in \mathcal{D}} P_\omega^i(\mathcal{T}_{m+3n+1}^\infty \in A \cap \mathcal{T}_1^m \in B)\mu_\omega(i)\right)$$

$$(44) \qquad = \left(E_P\left(\sum_{i \in \mathcal{D}} P_\omega(\mathcal{T}_{m+3n+1}^\infty \in A | \mathcal{T}_1^m \in B, X_0 = (0, i))\right.\right.$$

$$\left.\left. \times P_\omega^i(\mathcal{T}_1^m \in B)\mu_\omega(i)\right)\right)$$

$$\times [\mathbb{P}^\mu(\mathcal{T}_1^m \in B)]^{-1}$$

and

$$P_\omega(\mathcal{T}_{m+3n+1}^\infty \in A | \mathcal{T}_1^m \in B, X_0 = (0, i))$$

$$(45) \qquad = \sum_{j \in \mathcal{D}} P_\omega(\mathcal{T}_{m+3n=1}^\infty \in A | Y_{T_m} = j) \cdot P_\omega(Y_{T_m} = j | \mathcal{T}_1^m \in B, X_0 = (0, i))$$

$$= \sum_{j \in \mathcal{D}} P_{\theta^m\omega}^j(\mathcal{T}_{3n+1}^\infty \in A) \cdot P_\omega(Y_{T_m} = j | \mathcal{T}_1^m \in B, X_0 = (0, i)).$$

Next, $P_{\theta^m\omega}^j(\mathcal{T}_{3n+1}^\infty \in A) = \sum_{k \in \mathcal{D}} P_{\theta^m\omega}^j(Y_{T_{3n}} = k) \cdot P_{\theta^{m+3n}\omega}^k(\mathcal{T}_1^\infty \in A)$, and, using the notation of Lemma 3.3, we have $P$-a.s.,

$$|P_{\theta^m\omega}^j(Y_{T_{3n}} = k) - \pi_{m+3n}(j, k)| \leq c_m c_{m+1} \cdots c_{m+3n-1}, \qquad j, k \in \mathcal{D}.$$

Therefore, with probability 1, we have for all $j \in \mathcal{D}$,

$$|P^j_{\theta^m \omega}(\mathcal{T}^\infty_{3n+1} \in A) - P^{\pi_{\theta^{m+3n}\omega}}_{\theta^{m+3n}\omega}(\mathcal{T}^\infty_1 \in A)| \leq dc_m c_{m+1} \cdots c_{m+3n-1},$$

and hence, by using (45), $P$-a.s. for any $i \in \mathcal{D}$,

$$|P_\omega(\mathcal{T}^\infty_{m+3n+1} \in A | \mathcal{T}^m_1 \in B, X_0 = (0,i)) - P^{\pi_{\theta^{m+3n}\omega}}_{\theta^{m+3n}\omega}(\mathcal{T}^\infty_1 \in A)|$$
$$\leq dc_m \cdots c_{m+3n-1}.$$

By (44) this implies that

$$(46) \quad \left| \mathbb{P}^\mu(\mathcal{T}^\infty_{m+3n+1} \in A | \mathcal{T}^m_1 \in B) - \frac{E_P(P^{\pi_{\theta^{m+3n}\omega}}_{\theta^{m+3n}\omega}(\mathcal{T}^\infty_1 \in A)P^\mu_\omega(\mathcal{T}^m_1 \in B))}{\mathbb{P}^\mu(\mathcal{T}^m_1 \in B)} \right|$$
$$\leq \frac{d^2}{\mathbb{P}^\mu(\mathcal{T}^m_1 \in B)} \cdot E_P(c_m \cdots c_{m+3n-1} P^\mu_\omega(\mathcal{T}^m_1 \in B)).$$

We next show that the rightmost expression in (46) decays (as function of $n$) at least as $\varphi(n)$ [defined on (13)] with an exponentially small correction term, that is, there exist constants $K_1 > 0$ and $K_2 > 0$ such that for all $m \in \mathbb{N}$,

$$(47) \quad \frac{1}{\mathbb{P}^\mu(\mathcal{T}^m_1 \in B)} \cdot E_P(c_m \cdots c_{m+3n-1} P^\mu_\omega(\mathcal{T}^m_1 \in B))$$
$$< \varphi(n) + K_1 e^{-K_2 n}.$$

Note that $c_n = 1 - \max_{i \in \mathcal{D}} \min_{j \in \mathcal{D}} \eta_n(i,j)$ is a function of $(\omega_k)_{k \leq n}$ and therefore we cannot use the mixing property of the environment directly in order to separate $c_m \cdots c_{m+3n-1}$ and $P^\mu_\omega(\mathcal{T}^m_1 \in B)$. Rather, we shall approximate $c_i$ by some functions measurable with respect to the $\sigma$-algebra generated by a finite interval of the environment.

Define for $k \in \mathbb{N}$, $n \in \mathbb{Z}$ and $i,j \in \mathcal{D}$ the following functions of the environment $\omega$:

$$(48) \quad \begin{aligned} f_{n,k}(i) &= P_\omega(\chi_0 \leq n - k | X_0 = (n,i)), \\ \eta_{n,k}(i,j) &= P_\omega(Y_{T_1} = j, \chi_0 > n - k | X_0 = (n,i)), \\ c_{n,k} &= 1 - \max_{i \in \mathcal{D}} \min_{j \in \mathcal{D}} \eta_{n,k}(i,j), \\ \eta^-_n(i,j) &= P^i_{\theta^n \omega}(Y_{T_{-1}} = j). \end{aligned}$$

The numbers $c_{n,k}$ approximate $c_n$ and, in contrast to the latter, are measurable with respect to a $\sigma$-algebra generated by a finite interval of the environment [namely, w.r.t. $\sigma(\omega_i : n - k + 1 \leq i \leq n)$]. The quantity $\eta^-_n(i,j)$ is the probability that the walk starting at site $(n,i)$ reaches the layer $n - 1$ at point $(n-1,j)$. Note that $P$-a.s.,

$$c_n \leq c_{n,k} \leq c_n + \max_{i \in \mathcal{D}} f_{n,k}(i) \leq c_n + \|\eta^-_n \cdots \eta^-_{n-k+1}\|.$$

In particular this implies that for $k, m, n \in \mathbb{N}$,

$$
\begin{aligned}
&\prod_{i=m}^{m+n} c_{i,k} - \prod_{i=m}^{m+n} c_i \\
&\qquad \leq \left( c_{m,k} \prod_{i=m+1}^{m+n} c_{i,k} - c_{m,k} \prod_{i=m+1}^{m+n} c_i \right) + \left( c_{m,k} \prod_{i=m+1}^{m+n} c_i - c_m \prod_{i=m+1}^{m+n} c_i \right) \\
&\qquad \leq \left( \prod_{i=m+1}^{m+n} c_{i,k} - \prod_{i=m+1}^{m+n} c_i \right) + (c_{m,k} - c_m) \leq \cdots \\
&\qquad \leq \sum_{i=m}^{n+m} \| \eta_i^- \cdots \eta_{i-k+1}^- \|.
\end{aligned}
\tag{49}
$$

We need the following lemma.

LEMMA 5.4. *Let the conditions of Theorem 2.5 hold. Then, there exist constants $K_3 > 0$ and $K_4 > 0$ such that*

$$
E_P \left( \max_{i \in \mathcal{D}} P_\omega^i(\xi_k = -n \text{ for some } k \in \mathbb{N}) \right)
$$
$$
= E_P(\| \eta_0^- \cdots \eta_{-n+1}^- \|) \leq K_3 e^{-K_4 n}, \qquad n \in \mathbb{N}.
$$

PROOF. By Chebyshev's inequality and inequality (30), we have for all $n \in \mathbb{N}$,

$$
\begin{aligned}
E_P(\| \eta_0^- \cdots \eta_{-n+1}^- \|) &\leq E_P(\| M_{n,0} \|) \\
&\leq E_P(\| a_0 a_{-1} \cdots a_{-n+2} a_{-n+1}(\mathbf{1} + u_{-n}) \|),
\end{aligned}
$$

where the random vectors $M_{n,0}$ and $u_n$ are defined in (29) and (17), respectively.

Therefore, by the Cauchy–Schwarz inequality,

$$
E_P(\| \eta_0^- \cdots \eta_{-n+1}^- \|) \leq \sqrt{E_P(\| a_0 a_{-1} \cdots a_{-n+1} \|^2) E_P(\| \mathbf{1} + u_{-n} \|^2)}.
$$

Using the notation introduced in the proof of Lemma 3.6,

$$
\begin{aligned}
&E_P(\| u_{-n} \|^2) \\
&E_P(\| y_{-n+1} + \eta_{-n} y_{-n+2} + \eta_{-n} \eta_{-n+1} y_{-n+3} + \cdots + \eta_{-n} \cdots \eta_{-1} y_1 \|^2) \\
&\qquad \leq E_P([\| y_{-n+1} \| + \| y_{-n+2} \| + \cdots + \| y_1 \|]^2) \leq (n+1) E_P(\| y_0 \|^2),
\end{aligned}
\tag{50}
$$

where in the last step we use the fact that the sequence $(y_n)_{n \leq 0}$ is stationary in virtue of (21). In fact, similarly to (22), we have $y_n = b_n + a_n b_{n-1} + a_n a_{n-1} b_{n-2} + \cdots$.

Finally, using (22), Minkowski's inequality and the bound (10) for $\|b_n\|$, we get

$$(51) \quad [E_P(\|y_0\|^2)]^{1/2} \leq [E_P(\|b_0\|^2)]^{1/2} + \sum_{n=1}^{\infty} [E_P(\|a_0 \cdots a_{-n+1} b_{-n}\|^2)]^{1/2}$$

$$\leq 1/\epsilon \cdot \left(1 + \sum_{n=1}^{\infty} [E_P(\|a_0 \cdots a_{-n+1}\|^2)]^{1/2}\right).$$

The proof of the lemma is complete in view of assumption (ii) of Theorem 2.5.  $\square$

Coming back now to the rightmost expression of (46) we get, using the mixing and the translation invariance properties of the environment together with (49),

$$(52) \quad \frac{1}{\mathbb{P}^{\mu}(\mathcal{T}_1^m \in B)} \cdot E_P(c_m \cdots c_{m+3n-1} P_{\omega}^{\mu}(\mathcal{T}_1^m \in B))$$

$$\leq \frac{1}{\mathbb{P}^{\mu}(\mathcal{T}_1^m \in B)} \cdot E_P(c_{m+2n,n} \cdots c_{m+3n-1,n} P_{\omega}^{\mu}(\mathcal{T}_1^m \in B))$$

$$\leq E_P(c_{m+2n,n} \cdots c_{m+3n-1,n}) + \varphi(n)$$

$$\leq E_P(c_0 \cdots c_{n-1}) + n E_P(\|\eta_0^- \cdots \eta_{-n+1}^-\|) + \varphi(n),$$

where the mixing coefficients $\varphi(n)$ are defined in (13). This proves (47) in virtue of assumption (iv) of Theorem 2.5 and Lemma 5.4.

The last step of the proof is to show that there exist constants $K_5 > 0$ and $K_6 > 0$ such that for all $m \in \mathbb{N}$,

$$(53) \quad G(A, B) := \left| \frac{E_P(P_{\theta^{m+3n}\omega}^{\pi_{\theta^{m+3n}\omega}}(\mathcal{T}_1^{\infty} \in A) P_{\omega}^{\mu}(\mathcal{T}_1^m \in B))}{\mathbb{P}^{\mu}(\mathcal{T}_1^m \in B)} - \mathbb{P}^{\pi}(\mathcal{T}_1^{\infty} \in A) \right|$$

$$< K_5 e^{-K_6 n}.$$

This estimate combined together with (46) and (47) yields the claim of the lemma.

In order to prove the exponential upper bound for $G(A, B)$ we shall first approximate $P_{\theta^{m+3n}\omega}^i(\mathcal{T}_1^{\infty} \in A)$ by $P_{\theta^{m+3n}\omega}^i(\mathcal{T}_1^{\infty} \in A, \chi_{m+3n} > m+n) + P_{\theta^{m+3n}\omega}^i(\chi_{m+3n} \leq m+n)$ where the random variables $\chi_n$ are defined in (33). Note that both the summands in the approximation term are measurable with respect to $\sigma(\omega_i : i \geq m+n)$. Furthermore, for any $i \in \mathcal{D}$ and $n \in \mathbb{N}$ we have $P$-a.s.,

$$(54) \quad P_{\omega}^i(\mathcal{T}_1^{\infty} \in A) - P_{\omega}^i(\mathcal{T}_1^{\infty} \in A, \chi_0 > -n) \leq P_{\omega}^i(\chi_0 \leq -n)$$

$$= \eta_0^- \cdots \eta_{-n}^- \mathbf{1}(i).$$

Next, we shall approximate $\pi_{\theta^{m+3n}\omega}$ by

$$\pi_{\theta^{m+3n}\omega,n} := e_1\eta_{2n,n}\cdots\eta_{3n-1,n},$$

where the random matrices $\eta_{n,k}$ are defined in (48). For any $n \in \mathbb{N}$ we have

$$
\begin{aligned}
&\|\pi_\omega - \pi_{\omega,n}\| \\
(55) \quad &\leq \|\pi_\omega - e_1\eta_{-n}\cdots\eta_{-1}\| + \|e_1\eta_{-n}\cdots\eta_{-1} - e_1\eta_{-n,n}\cdots\eta_{-1,n}\| \\
&\leq c_{-n}\cdots c_{-1} + \sum_{i=1}^n \|\eta_{-i} - \eta_{-i,n}\| \\
&\leq c_{-n}\cdots c_{-1} + \sum_{i=1}^n \|\eta_{-i}^-\cdots\eta_{-i-n+1}^-\|,
\end{aligned}
$$

where in the last but one step we use Lemma 3.3 in order to bound the first term and the iteration argument similar to the one leading to (49) in order to bound the second one.

It follows from (54) and (55) that

$$
\begin{aligned}
0 &\leq P_\omega^\pi(\mathcal{T}_1^\infty \in A) - \sum_{i=1}^d \pi_{\omega,n}(i)P_\omega^i(\mathcal{T}_1^\infty \in A, \chi_0 \geq m+2n) \\
(56) \quad &= \sum_{i=1}^d (\pi_\omega(i) - \pi_{\omega,n}(i))P_\omega^i(\mathcal{T}_1^\infty \in A) \\
&\quad + \sum_{i=1}^d \pi_{\omega,n}(i)(P_\omega^i(\mathcal{T}_1^\infty \in A) - P_\omega^i(\mathcal{T}_1^\infty \in A, \chi_0 \geq m+2n)) \\
&\leq d\left[c_{-n}\cdots c_{-1} + \sum_{i=1}^n \|\eta_{-i} - \eta_{-i,n}\| + \|\eta_0^-\cdots\eta_{-n}^-\|\right].
\end{aligned}
$$

Therefore, using (47) and Lemma 30, we obtain that the following inequality holds with suitable constants $K_7, K_8$ and $K_9$:

$$
\begin{aligned}
&G(A,B) \\
&\quad \leq \left|\frac{E_P(\sum_{i=1}^d \pi_{\theta^{m+3n}\omega,n}(i)P_{\theta^{m+3n}\omega}^i(\mathcal{T}_1^\infty \in A)P_\omega^\mu(\mathcal{T}_1^m \in B))}{\mathbb{P}^\mu(\mathcal{T}_1^m \in B)} - \mathbb{P}^\pi(\mathcal{T}_1^\infty \in A)\right| \\
&\qquad + K_7\varphi(n) + K_8 e^{-K_9 n} \\
&\quad \leq \mathbb{P}^\pi(\mathcal{T}_1^\infty \in A) - E_P\left(\sum_{i=1}^d \pi_{\omega,n}(i)P_\omega^i(\mathcal{T}_1^\infty \in A, \chi_0 \geq -n)\right) \\
&\qquad + (K_7+1)\varphi(n) + K_8 e^{-K_9 n},
\end{aligned}
$$

which proves (53) in view of (56). The proof of Proposition 5.3 is thus complete.  □

The next lemma shows that $\mathbb{E}^\pi(\tau_1^2) < \infty$ under the conditions of Theorem 2.5.

LEMMA 5.5.    *Let the conditions of Theorem 2.5 hold. Then $\mathbb{E}^\pi(\tau_1^2) < \infty$.*

PROOF.    Consider the RWRE $(X_k)_{k \geq 0}$ with $X_0 = (n, i)$, where $n \leq 0$ and $i \in \mathcal{D}$. Taking the square in both sides of identity (18) yields

$$\tau_1^2 = 1 + \sum_{k \in \{-1,0,1\}} \sum_{j \in \mathcal{D}} \mathbf{I}_{\{\xi_1 = n+k, Y_1 = j\}} (\tilde\tau_1 + 2\tilde\tau_1).$$

Recall the vectors $u_n$ defined in (17). Taking in the last identity quenched expectations and using the notation $w_n^i = E_\omega(T_1^2 | X_0 = (n, i))$ and $w_n = (w_n^1, \ldots, w_n^d)$, we obtain

$$
\begin{aligned}
(57) \qquad w_n &= \mathbf{1} + q_n w_{n-1} + r_n w_n + p_n w_{n+1} + 2q_n u_{n-1} + 2r_n u_n + 2p_n u_{n+1} \\
&= q_n w_{n-1} + r_n w_n + p_n w_{n+1} + 2u_n - \mathbf{1}.
\end{aligned}
$$

Let $z_1 := w_0$ and for $n \leq 0$, $z_n := w_{n-1} - \eta_{n-1} w_n$ using the convention that $\infty - \infty = \infty$. Substituting $w_{n-1} = z_n + \eta_{n-1} w_n$ in (57) we obtain, similarly to (21) and (22),

$$z_{n+1} = \gamma_n(2u_n - \mathbf{1}) + a_n z_n, \qquad n \geq 0,$$

and

$$w_0 = z_1 = \gamma_0(2u_0 - \mathbf{1}) + \sum_{n=1}^\infty a_0 a_{-1} \cdots a_{-n+1} \gamma_{-n}(2u_{-n} - \mathbf{1}).$$

To complete the proof it suffices to show that there are positive constants $K_5 > 0$ and $K_6 > 0$ such that for any $n \geq 0$,

$$(58) \qquad E_P(\|a_0 a_{-1} \cdots a_{-n+1} \gamma_{-n} u_{-n}\|) \leq K_5 e^{-K_6 n}.$$

Toward this end we first observe that

$$
\begin{aligned}
&E_P(\|a_0 a_{-1} \cdots a_{-n+1} \gamma_{-n} u_n\|) \\
&\qquad \leq E_P(\|a_0 a_{-1} \cdots a_{-n+1}\| \cdot \|u_{-n}\|) \\
&\qquad \leq [E_P(\|a_0 a_{-1} \cdots a_{-n+1}\|^2)]^{1/2} \cdot [E_P(\|u_{-n}\|^2)]^{1/2}.
\end{aligned}
$$

The estimate (58) with suitable constants $K_5$ and $K_6$ follows now from inequalities (50) and (51) together with assumption (ii) of Theorem 2.5. This completes the proof of the lemma.  □

To conclude the proof of Theorem 2.5 it remains to show that $\sigma^2$ defined in (43) is strictly positive for $Z_n = \tau_n - v_P^{-1}$. To this end, we shall use a representation of $\sigma^2$ as a sum of a "quenched" and an "annealed" contributions and will show that the "annealed" term must be strictly positive. Similar methods of estimating the limiting variance have been employed, for instance, in [7] (see Theorem 4.3 and especially (26) there) and [24] (see the proof of Lemma 3.4). In the case of the strip model, as it is first observed in [15], it is convenient to think that the "annealed" fluctuations are due both to the randomness of the environment $\omega$ and to the randomness of the sequence $(Y_{T_n})_{n \geq 1}$.

Define for $n \geq 0$ and $i, j \in \mathcal{D}$,

$$u_n^{i,j}(\omega) = E_\omega(\tau_1 | X_0 = (n, i), Y_{\tau_1} = j) - v_P^{-1},$$

and observe that, since the random variables $\tau_n$ are independent under the conditional measure $P_\omega(\cdot | (Y_{T_n})_{n \geq 0})$,

$$\sigma^2 = \mathbb{E}^\pi((\tau_1 - v_P^{-1})^2) + 2 \sum_{n=1}^\infty \mathbb{E}^\pi(u_1^{Y_0, Y_{T_1}}(\omega) u_n^{Y_{T_n}, Y_{T_{n+1}}}(\omega))$$

$$= \mathbb{E}^\pi(\tau_1^2) - \mathbb{E}^\pi[(E_\omega(\tau_1 | X_0 = (n, i), Y_{\tau_1} = j))^2]$$

$$+ \mathbb{E}^\pi[(u_1^{Y_0, Y_{T_1}}(\omega))^2] + 2 \sum_{n=2}^\infty \mathbb{E}^\pi(u_1^{Y_0, Y_{T_1}}(\omega) u_n^{Y_{T_n}, Y_{T_{n+1}}}(\omega)).$$

It follows that $\sigma^2 > 0$ because Lemma 20.3 in [4], page 172, implies that

$$\mathbb{E}^\pi[(u_1^{Y_0, Y_{T_1}}(\omega))^2] + 2 \sum_{n=2}^\infty \mathbb{E}^\pi(u_1^{Y_0, Y_{T_1}}(\omega) u_n^{Y_{T_n}, Y_{T_{n+1}}}(\omega))$$

$$= \lim_{n \to \infty} \frac{1}{n} \mathbb{E}^\pi \left[ \left( \sum_{i=1}^n u_n^{Y_{T_n}, Y_{T_{n+1}}}(\omega) \right)^2 \right] \geq 0,$$

while Jensen's inequality together with the first inequality in (5) imply that

$$\mathbb{E}^\pi(\tau_1^2) - \mathbb{E}^\pi[(E_\omega(\tau_1 | X_0 = (n, i), Y_{\tau_1} = j))^2] > 0.$$

The proof of Theorem 2.5 is complete.

**Acknowledgments.** During the preparation of this paper, I became aware of work by Ilya Goldsheid on a similar subject. I am greatly indebted to Professor Goldsheid for the interesting and stimulating discussions on this topic and for showing me a copy of his paper [15] prior to its publication. I am very grateful to the referees for the careful reading of this paper and many helpful remarks and suggestions. I wish also to thank Jon Peterson for drawing my attention to an inconsistency of a definition in a preliminary version of the paper.

DEPARTMENT OF MATHEMATICS
IOWA STATE UNIVERSITY
AMES, IOWA 50011
USA
E-MAIL: roiterst@iastate.edu